\documentclass[11pt,reqno]{amsart}
\usepackage[T1]{fontenc}
\usepackage[latin9]{inputenc}
\usepackage{geometry}
\usepackage{verbatim}
\usepackage{mathrsfs}
\usepackage{amsthm}
\usepackage{amstext}
\usepackage{amssymb}
\usepackage{graphicx}
\usepackage{esint}
\usepackage{color}
\usepackage{amsfonts}
\usepackage{amsmath} 
\usepackage{bm}           
\usepackage{mathabx}
\usepackage{amsxtra}
\usepackage{stmaryrd}
\usepackage{old-arrows}
\usepackage{float,epsf}
\usepackage[english]{babel}
\usepackage{enumerate}
\usepackage{tikz}
\usepackage{srcltx}
\usepackage{verbatim}
\usepackage[multiple]{footmisc}

\makeatletter

\newcommand{\R}{\mathbb{R}}

\newcommand{\D}{{\rm D}}

\newcommand{\M}{\mathcal{M}}

 \newcommand{\xx}{\mathbf{x}}

\newcommand{\del}{{\partial}} 
\renewcommand{\d}{\delta}

\newcommand{\s}{\sigma} 
\renewcommand{\k}{\kappa}

\newcommand{\V}{{\mathbf{V}}}

\theoremstyle{plain}
\newtheorem{theorem}{Theorem}[section]

\newtheorem{problem}[theorem]{Problem}

\newcommand{\nnu}{{\boldsymbol \nu}}

\newcommand{\vv}{{\boldsymbol v}}

\def\intave#1{\int_{#1}\hbox{\llap{$\raise2.3pt\hbox{\vrule
height.9pt width7pt}\phantom{\scriptstyle{#1}}\mkern-2mu$}}}

\def\intav#1{\mathchoice
          {\mathop{\vrule width 9pt height 3 pt depth -2.6pt
                  \kern -9pt \intop}\nolimits_{\kern -6pt#1}}%
          {\mathop{\vrule width 5pt height 3 pt depth -2.6pt
                  \kern -6pt \intop}\nolimits_{#1}}%
          {\mathop{\vrule width 5pt height 3 pt depth -2.6pt
                  \kern -6pt \intop}\nolimits_{#1}}%
          {\mathop{\vrule width 5pt height 3 pt depth -2.6pt
                  \kern -6pt \intop}\nolimits_{#1}}}
\def\intav#1{\vint_{#1}}

\allowdisplaybreaks[1]
\allowdisplaybreaks[1]

\newcommand{\bR}{{\mathbb R}}

\newcommand{\sS}{\mathcal S}

\newcommand{\divg}{ \mbox{div\,}}

\newcommand{\vphi} {\varphi}

\newcommand{\x}{{\bf x}}

\renewcommand{\d}{\partial}
\newcommand{\G}{\Gamma}
\renewcommand{\k}{\kappa}
\renewcommand{\r}{\rho}

\newcommand{\bn}{{\boldsymbol\nu}}
\newcommand{\mR}{\mathbb{R}}
\newcommand{\CS}{\mathcal{S}}

\newcommand{\pSi}{\varphi}
\newcommand{\PtIncW}{{P_0}}
\newcommand{\PtUpL}{{P_1}}
\newcommand{\PtLwL}{{P_2}}
\newcommand{\PtLwR}{{P_3}}
\newcommand{\PtUpR}{{P_4}}

\newcommand{\xxi}{{\boldsymbol\xi}}
\newcommand \Sonic{\Gamma_{\rm sonic}}
\newcommand \Shock{\Gamma_{\rm shock}}
\newcommand \Wedge{\Gamma_{\rm wedge}}
\newcommand \Symm{\Gamma_{\rm sym}}
\newcommand{\cxi}{{\xi_1}}

\newcommand \shock{\Gamma_{\rm shock}}

\newcommand \ol{\overline}

\newcommand{\e}{{\varepsilon}}

\usepackage{babel}

\usepackage{babel}

\makeatother

\usepackage{babel}

\numberwithin{equation}{section}

\begin{document}
\title[Partial Differential Equations of Mixed Type]{Partial Differential Equations of Mixed Type\\ --- Analysis and Applications}
\author{Gui-Qiang G. Chen}
\address{Gui-Qiang G. Chen: Oxford Centre for Nonlinear Partial Differential Equations,
Mathematical Institute, University of Oxford, Oxford, OX2 6GG, UK.
 $\quad$ Email: chengq@maths.ox.ac.uk}
\maketitle

Partial differential equations (PDEs) are at the heart of many
mathematical and
scientific advances.
While great progress has been made on the theory of PDEs of standard types
during the last eight decades,
the analysis of nonlinear PDEs of mixed type is still in its infancy.
The aim of this expository paper is to show, through several longstanding fundamental problems
in fluid mechanics, differential geometry, and other areas,
that many nonlinear PDEs arising in these areas
are no longer of standard types, but
lie at the boundaries of the classification of PDEs or, indeed, go beyond the classification to be
of mixed type.
Some interrelated connections, historical perspectives, recent developments,
and current trends in the analysis of nonlinear PDEs of mixed type are also presented.

\section{Linear Partial Differential Equations of Mixed Type}

\begin{figure}
\vspace{-8pt}
\begin{minipage}{0.45\textwidth}
\centering
\includegraphics[height=1.39in,width=1.20in]{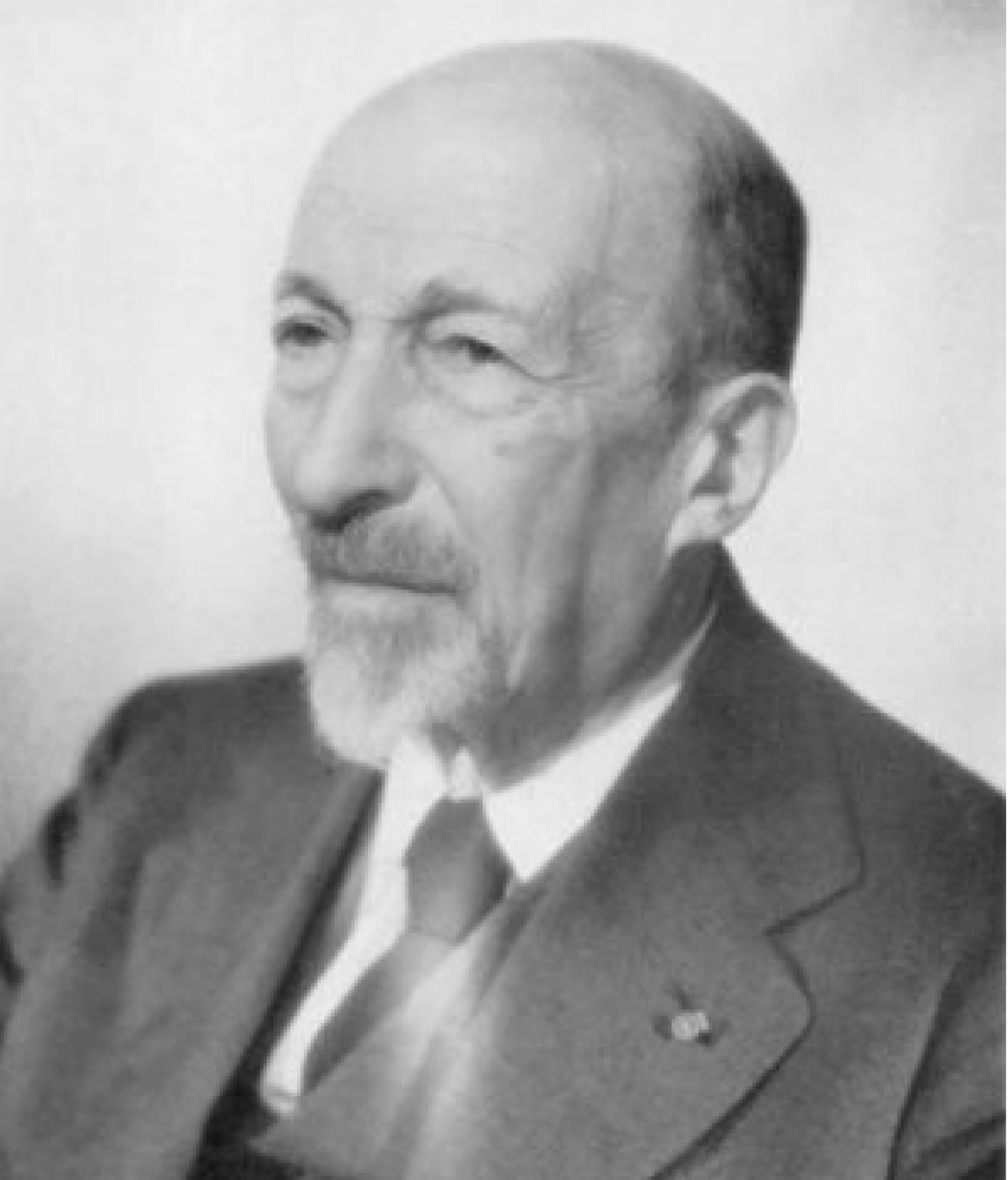}
\caption[]{
\scriptsize
\,\,\,\,\,\,\,\,Jacques Salomon Hadamard (8 December 1865 -- 17 October 1963)
first introduced the classification of PDEs in \cite{Had}}\label{figure-1}
\end{minipage}
\begin{minipage}{0.51\textwidth}
\centering
\includegraphics[height=1.35in,width=2.65in]{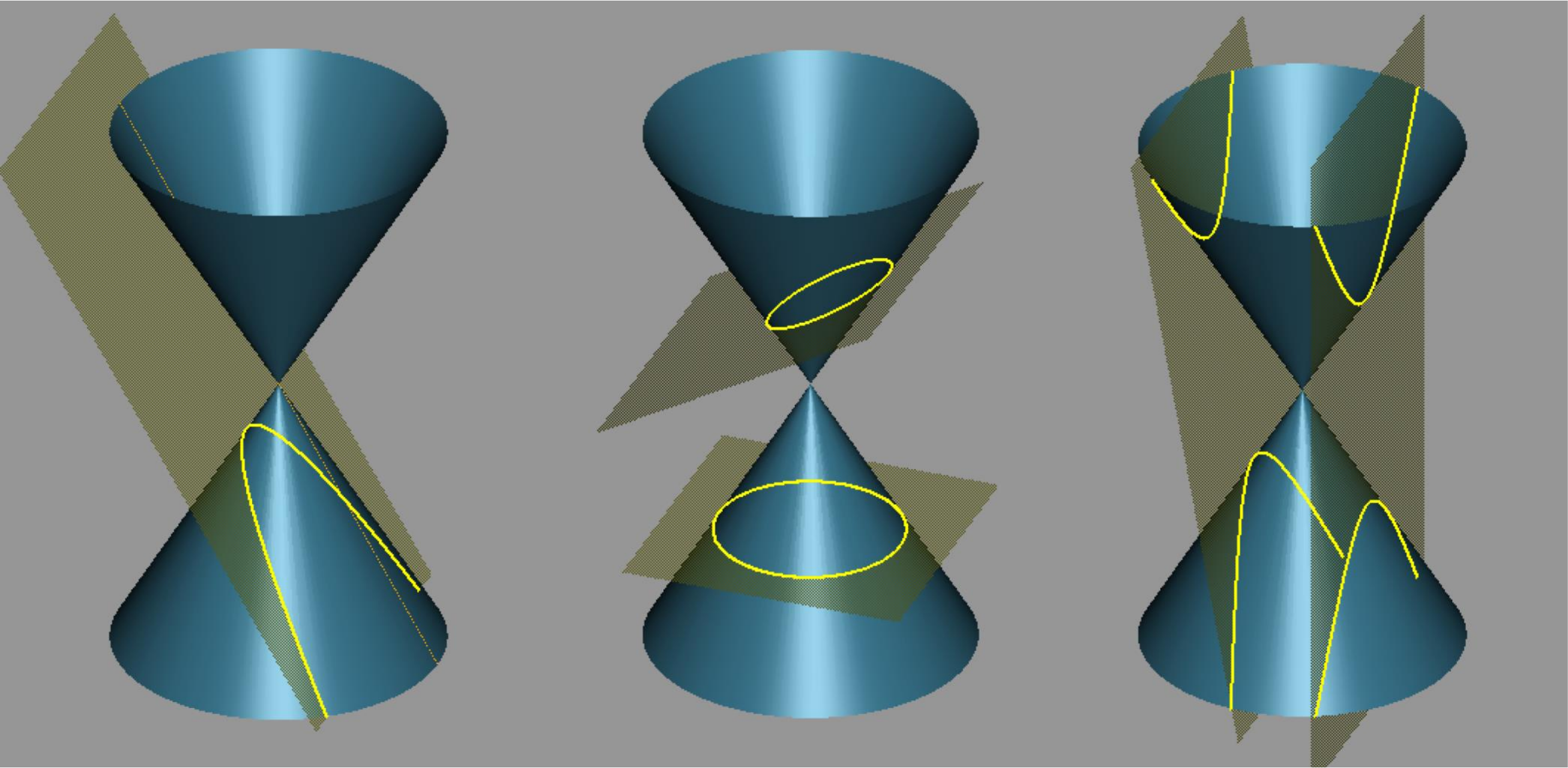}
\caption[]{
\scriptsize
\,\, Types of conic sections: Parabola, ellipse, and hyperbola}\label{figure-2}
\end{minipage}
\end{figure}

Three of the basic types of PDEs are
the {\it elliptic}, {\it hyperbolic}, and {\it parabolic} types, following
the classification
introduced by Jacques  Salomon  Hadamard in 1923
(see Fig. \ref{figure-1}).

The prototype of second-order elliptic equations is the {\it Laplace equation}:
\begin{equation}\label{1.b}
\Delta_\x u:=\sum_{j=1}^n  \partial_{x_jx_j}u=0 \qquad\quad\mbox{for $\x=(x_1, \dots, x_n)\in \mathbb{R}^{n}$},
\end{equation}
often describing physical equilibrium states, whose solutions are also called harmonic functions
or potential functions, where $\partial_{x_jx_j}$ is the second-order partial derivative in the $x_j$-variable, $j=1,\dots, n$.
The simplest representative of  hyperbolic equations is the {\it wave equation}:
\begin{equation}\label{1.a}
\partial_{tt}u-\Delta_\x u=0
\qquad\quad \mbox{for $(t, \x)\in \mathbb{R}^{n+1}$},
\end{equation}
which governs the propagation of linear waves
(such as acoustic waves and electromagnetic waves),
while the prototype of second-order parabolic equations is the {\it heat equation}:
\begin{equation}\label{1.c}
\partial_t u-\Delta_\x u=0  \qquad\quad \mbox{for $(t, \x)\in \mathbb{R}^{n+1}$},
\end{equation}
which often describes the dynamics of temperature and diffusion/stochastic processes.

At first glance, the forms of the Laplace/heat equations and the wave equation look quite similar.
In particular, any steady solution of the wave/heat equations is a solution of the Laplace equation,
while a solution of the Laplace equation often determines an asymptotic state of the time-dependent
solutions of the wave/heat equations.
However, the properties of solutions of the Laplace/heat equations and the wave equation
are significantly different.
One difference is the {\it infinite versus finite speed of propagation} of the solution, while another
is the {\it gain versus loss of regularity} of the solution; see \cite{Evans,Had} and the references cited therein.
Since the solutions of elliptic/parabolic PDEs share many common features, we will focus mainly on
PDEs of mixed elliptic-hyperbolic type
from now on.

\smallskip
The distinction between the elliptic and hyperbolic types can be seen more clearly
from the classification of two-dimensional ($2$-D) constant-coefficient second-order PDEs:
\begin{equation}\label{1.1}
a_{11} \partial_{x_1x_1}u+2a_{12}\partial_{x_1x_2}u +a_{22}\partial_{x_2x_2}u=f(\x) \qquad\quad\,\, \mbox{for $\,\x=(x_1,x_2)\in \R^2$}.
\end{equation}
Let $\lambda_1\le \lambda_2$ be the two constant eigenvalues of the $2\times 2$ symmetric coefficient
matrix $(a_{ij})_{2\times 2}$.
Eq. \eqref{1.1} is classified to be {\it elliptic} if
\begin{equation}\label{1.2}
\det(a_{ij})>0 \iff \lambda_1\lambda_2>0  \iff a_{12}^2-a_{11}a_{22}<0,
\end{equation}
while Eq. \eqref{1.1} is classified to be {\it hyperbolic} if
\begin{equation}\label{1.3}
\det(a_{ij})<0 \iff \lambda_1\lambda_2<0  \iff a_{12}^2-a_{11}a_{22}>0.
\end{equation}
Notice that the left-hand side of Eq. \eqref{1.1} is analogous to the quadratic (homogeneous) form:
$$
 a_{11} \xi_1^2+ a_{12}\xi_1\xi_2+a_{22}\xi_2^2
$$
for conic sections.
Thus, the classification of Eq. \eqref{1.1} is consistent with the classification of
conic sections and quadratic forms
in algebraic geometry,
based on the sign of the discriminant: $a_{12}^2-a_{11}a_{22}$.
The corresponding quadratic curves are
ellipses (incl. circles), hyperbolas, and parabolas (see Fig. \ref{figure-2}).

This classification can also be seen through taking the Fourier transform
on both sides of Eq. \eqref{1.1}:
\begin{equation}\label{1.4}
(a_{11}\xi_1^2+2a_{12}\xi_1\xi_2 +a_{22}\xi_2^2)\, \hat{u}(\xxi)= -\hat{f}(\xxi)
\qquad\quad\,\, \mbox{for $\,\xxi=(\xi_1, \xi_2)\in \mathbb{R}^2$},
\end{equation}
where
$
\hat{w}(\xxi)=\frac{1}{2\pi}\int_{-\infty}^\infty w(\x)e^{-i \x\cdot\xxi}\, {\rm d}\x
$
is the Fourier transform
of a function $w(\x)$, such as $u(\x)$ and $f(\x)$ for \eqref{1.4}.
When Eq. \eqref{1.1} is {\it elliptic}, the Fourier transform $\hat{u}(\xxi)$ of
the solution $u(\x)$ gains two orders of decay
for the high Fourier frequencies ({\it i.e.}, $|\xxi|\gg 1$)
so that the solution gains the regularity of two orders from $f(\x)$.
When Eq. \eqref{1.1} is {\it hyperbolic},
$\hat{u}(\xxi)$ fails to gain two orders of decay
for the high Fourier frequencies along the two characteristic directions
in which
$a_{11}\xi_1^2+2a_{12}\xi_1\xi_2 +a_{22}\xi_2^2=0$,
even though it still gains two orders of decay for the high Fourier frequencies
away from these two
directions.

For the classification above, a general homogeneous constant-coefficient second-order
PDE ({\it i.e.}, $f(\x)=0$)
with \eqref{1.2} or \eqref{1.3} can be transformed
into the Laplace equation \eqref{1.b} with $n=2$,
or the wave equation \eqref{1.a} with $n=1$, via the corresponding
coordinate transformations, respectively.
This exposes the beauty of the classification theory introduced first in Hadamard \cite{Had}.

\medskip
On the other hand, for general variable-coefficient second-order PDEs:
\begin{equation}\label{1.5}
a_{11}(\mathbf{x})\partial_{x_1x_1}u +2 a_{12}(\mathbf{x})\partial_{x_1x_2}u +a_{22}(\mathbf{x})\partial_{x_2x_2}u=f(\x),
\end{equation}
the situation is different.
The classification depends upon the signature of the eigenvalues $\lambda_j(\mathbf{x}), j=1,2$,
of the coefficient matrix $(a_{ij}(\mathbf{x}))$.
In general,
$
\lambda_1(\mathbf{x})\lambda_2(\mathbf{x})
$
may change its sign as a function of $\x$, which leads to the  {\it mixed elliptic-hyperbolic type} of \eqref{1.5}:
Eq. \eqref{1.5} is {\it elliptic} when $\lambda_1(\mathbf{x})\lambda_2(\mathbf{x})>0$ and
{\it hyperbolic} when  $\lambda_1(\mathbf{x})\lambda_2(\mathbf{x})<0$ with transition boundary/region where
$\lambda_1(\mathbf{x})\lambda_2(\mathbf{x})=0$.

\smallskip
Three of the classical prototypes for linear PDEs of mixed elliptic-hyperbolic type are:

\smallskip
(i)
{\it The Lavrentyev-Bitsadze equation}:
$\qquad\partial_{x_1x_1}u +{\rm sign} (x_1)\partial_{x_2x_2}u=0$.

\noindent
This equation exhibits
a jump transition at $x_1=0$.
It becomes the Laplace equation \eqref{1.b} in the half-plane $x_1>0$ and
the wave equation \eqref{1.a} in the half-plane $x_1<0$, and changes its type from {\it elliptic}
to {\it hyperbolic} via the jump-discontinuous coefficient ${\rm sign} (x_1)$.

\smallskip
(ii) {\it The Tricomi equation}:
$\qquad \partial_{x_1x_1}u + x_1\partial_{x_2x_2}u =0$.

\noindent
This equation is of hyperbolic degeneracy at $x_1=0$.
It is {\it elliptic} in the half-plane $x_1>0$ and {\it hyperbolic} in the half-plane $x_1<0$, and
changes its type from {\it elliptic} to {\it hyperbolic} through the degenerate line $x_1=0$.
This equation is of hyperbolic degeneracy in the domain $x_1\le 0$, in which the
two characteristic families coincide {\it perpendicularly} to the line $x_1=0$. Its degeneracy is determined
by the classical elliptic or hyperbolic {\it Euler-Poisson-Darboux equation}\footnote{Hadamard, J.:
{\it La Th\'{e}orie des \'{E}quations aux D\'{e}riv\'{e}es Partielles},
in French, \'{E}ditions Scientifiques: Peking; Gauthier-Villars \'{E}diteur:
Paris, 1964.}:
\begin{equation}\label{1.8}
\partial_{\tau\tau}u \pm \partial_{x_2x_2}u +\frac{\beta}{\tau}\partial_\tau u =0
\end{equation}
with $\beta=\frac{1}{3}$ for $\tau=\frac{2}{3}|x_1|^{\frac{3}{2}}$,
and signs ``$\pm$''
corresponding to the
half-planes $\pm x_1>0$ for $\x$ to lie in.

\smallskip
(iii) {\it The Keldysh equation}:
$\qquad x_1 \partial_{x_1x_1}u + \partial_{x_2x_2}u =0$.

\noindent
This equation is of parabolic degeneracy at $x_1=0$.
It is {\it elliptic} in the half-plane $x_1>0$ and {\it hyperbolic} in the half-plane $x_1<0$, and
changes its type from elliptic to hyperbolic through the degenerate line $x_1=0$.
This equation is of parabolic degeneracy in the domain $x_1\le 0$, in which the
two characteristic families are quadratic parabolas lying in the half-plane $x_1<0$, and tangential at contact points
to the degenerate line $x_1=0$.
Its degeneracy is also determined
by the classical elliptic or hyperbolic {\it Euler-Poisson-Darboux equation}
\eqref{1.8}
with $\beta=-\frac{1}{4}$ for $\tau=\frac{1}{2}|x_1|^{\frac{1}{2}}$.

\smallskip
For such a linear PDE, the transition boundary
({\it i.e.}, the boundary between the elliptic and hyperbolic domains) is {\it a priori}
known. Thus, one traditional approach is to regard such a PDE as a degenerate elliptic
or hyperbolic PDE in the corresponding domain, and then to analyze
the solution behavior of these degenerate PDEs separately in the elliptic and hyperbolic domains
with degeneracy on the transition boundary, determined by the Euler-Poisson-Darboux type equations,
say \eqref{1.8}.
Another successful approach for dealing with such a PDE is the fundamental
solution approach.
That is, we first understand the behavior of the fundamental solution
of the mixed-type PDE,
especially its singularity,
from which analytical/geometric properties of solutions can then be revealed,
since the fundamental solution is a generator of all the solutions of the linear PDE.
Great effort and progress have been made in the analysis of linear PDEs of mixed type
by many leading mathematicians
since the early 20th century ({\it cf.} \cite{Bet,CF18,Had,Mo85}
and the references cited therein).
Still there are many important problems for linear PDEs of mixed type,
which require further understanding.

\smallskip
In the next sections, we show, through several longstanding fundamental
problems in  fluid mechanics, differential geometry, and other areas,
that many nonlinear PDEs arising in mathematics
and science are no longer of standard type, but are of mixed type.
In contrast to the linear case, the transition boundary for a nonlinear PDE of mixed type
is {\it a priori} unknown in general, and the nonlinearity often generates additional singularities.
Thus,
many classical methods and techniques
for linear PDEs no longer work for
nonlinear PDEs of mixed type directly.
The lack of effective unified approaches
is one of the main obstacles for tackling
the elliptic/hyperbolic phases together for
nonlinear PDEs of mixed type.
During the last eight decades, the PDE research community has been largely partitioned
by the approaches taken to
the analysis of different classes of PDEs (elliptic/hyperbolic/parabolic).
However, advances in the analysis of nonlinear PDEs over the
last several decades have made
it increasingly clear that many difficult
questions faced by the community are at the boundaries of this classification
or, indeed, go beyond this
classification.
In particular, many important nonlinear PDEs that arise in longstanding fundamental
problems across diverse areas
are of mixed type.
As we will show in \S 2--\S 4 below,
these problems include steady transonic flow problems and
shock reflection/diffraction problems in gas dynamics,
high-speed flow,
and related areas
({\it cf}. \cite{BD,Bers,CF18,CFr,Da,GM,Mo85,vN1,vonN}),
and isometric embedding problems with optimal target dimensions
and assigned regularity/curvatures
in elasticity, geometric analysis, materials science,
and other areas
({\it cf}.
\cite{CSW10,HH}).
The solution to these problems
will
advance our understanding of shock reflection/diffraction phenomena,
transonic flows, properties/classifications of elastic/biological surfaces/bodies/manifolds,
and other scientific issues,
and lead to significant developments of these areas and related mathematics.
To achieve this, a deep understanding of the underlying nonlinear PDEs of mixed type
(for instance, the solvability, the properties of solutions, {\it etc}.) is key.

\section{Nonlinear PDEs of Mixed Type and Steady Transonic Flow Problems in Fluid Mechanics}
In many applications, fluid flows are often regarded as time-independent,
which is the case for some longstanding fundamental problems,
such as the transonic flows past
multi-dimensional (M-D) obstacles (wedges/conic bodies, airfoils, {\it etc.}) or de Laval nozzles; see Figs. 3--4.
Furthermore, steady-state solutions are often global attractors as time-asymptotic equilibrium states
and serve as building blocks for constructing
time-dependent solutions ({\it cf.} \cite{CF18,CFr,Da,GM}).
The underlying nonlinear PDEs governing these fluid flows are generically
of mixed type.

Our first example is
steady potential fluid flows governed
by the steady Euler equations of conservation law of mass and Bernoulli's law:
\begin{equation}\label{potential-eq}
\divg(\rho\, \nabla\varphi)=0, \quad
{1\over 2}|\nabla\varphi|^2+{1\over \gamma-1}\rho^{\gamma-1}=\frac{B_0}{\gamma-1}
\qquad\quad\,\, \mbox{for ${\bf x}\in \bR^n$},
\end{equation}
after scaling, where $\rho$ is the density, $\varphi$ the velocity potential ({\it i.e.}, $\vv=\nabla\varphi$ is the velocity),
$\gamma>1$ the adiabatic exponent for
the ideal gas,  $B_0/(\gamma-1)$ the Bernoulli constant,
and $\nabla$ the gradient in
$\xx$.
System \eqref{potential-eq}, along with its time-dependent version (see \eqref{potential-t-eq} below),
is among the first PDEs to be written down by Euler ({\it cf.} Fig. \ref{figure-6a}), and
has been used widely in aerodynamics and other areas when the
vorticity waves are weak in the fluid flow under consideration
({\it cf.} \cite{Bers,CF18,CFr,Da,GM}).
System \eqref{potential-eq} for the steady
velocity potential $\varphi$
can be rewritten as
\begin{equation}\label{steady-potential-eq}
\divg\big(\rho_B(|\nabla\varphi|)\nabla\varphi\big)=0
\qquad\quad\quad\mbox{with $\rho_B(q)=\big(B_0-(\gamma-1)q^2/2\big)^{{1}/(\gamma-1)}$}.
\end{equation}
Eq.\,\eqref{steady-potential-eq} is a {\it nonlinear conservation law of
mixed elliptic-hyperbolic type} -- It is
\begin{itemize}
\item  strictly {\it elliptic} (subsonic)
if
$|\nabla\varphi|<c_*:=\sqrt{2B_0/(\gamma+1)}$,

\smallskip
\item strictly {\it hyperbolic} (supersonic) if
$|\nabla\varphi|>c_*$.
\end{itemize}
The transition boundary
is $|\nabla\varphi|=c_*$ (sonic),
a degenerate set of \eqref{steady-potential-eq},
which is {\it a priori} unknown,
since it is determined by the solution itself.

Similarly, the time-independent full Euler flows are governed by
the steady Euler equations:
\begin{equation}\label{steady-euler-syst}
\divg(\rho \vv)=0,\quad \divg(\rho \vv\otimes \vv) +\nabla p=0,
\quad \divg\big(\rho \vv(E+\frac{p}{\rho})\big)=0,
\end{equation}
where $p$ is the pressure, $\vv$ the velocity,
and $E=\frac{1}{2}|\vv|^2+e$ the energy with $e=\frac{p}{(\gamma-1)\rho}$ as
the internal energy.
System \eqref{steady-euler-syst} is a {\it system of conservation laws of mixed-composite hyperbolic-elliptic type} -- It is
\begin{itemize}
\item strictly {\it hyperbolic} when $|\vv|>c$ (supersonic),

\smallskip
\item {\it mixed-composite elliptic-hyperbolic} (two of them are elliptic and the others are hyperbolic) when
$|\vv|< c$ (subsonic),
\end{itemize}
where $c=\sqrt{\gamma p/\rho}$ is the sonic speed.
The transition boundary
between the supersonic/subsonic phase
is
$|\vv|=c$,
a degenerate set of the solution of system \eqref{steady-euler-syst},
which is {\it a priori} unknown.
\begin{figure}
\begin{minipage}{0.485\textwidth}
\centering
\includegraphics[height=0.75in,width=1.55in]{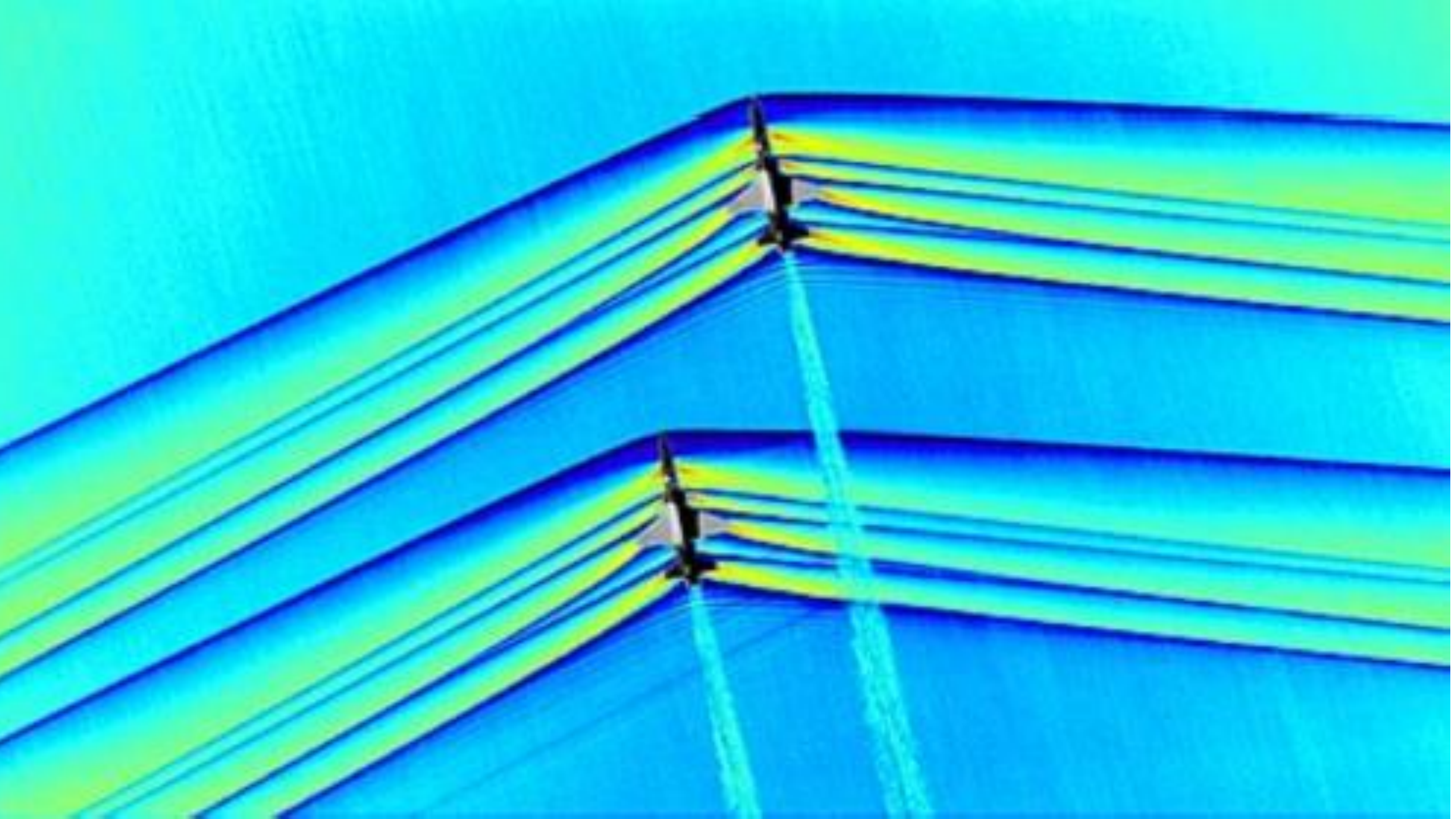}
\caption{\scriptsize \,\,
NASA's first Schlieren photo of shock waves interacting between two aircrafts
taken in
March 2019} \label{figure-3a}
\end{minipage}
\begin{minipage}{0.485\textwidth}
\centering
\includegraphics[height=0.70in,width=1.9in]{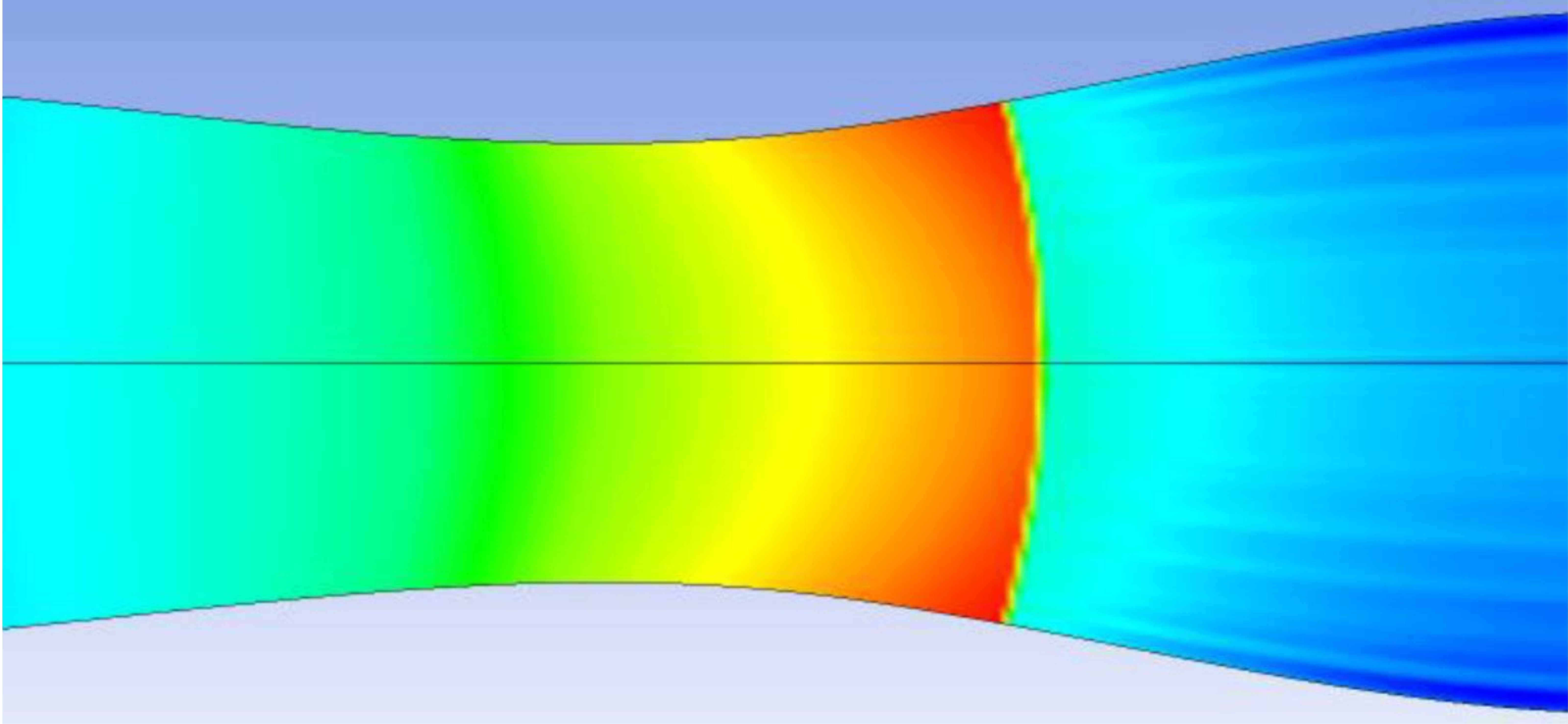}
\caption{\scriptsize Diagram of a de Laval nozzle for the approximate flow velocity}
\label{figure-3b}
\end{minipage}
\end{figure}

\begin{figure}
\begin{minipage}{0.485\textwidth}
\vspace{-10pt}
\centering
\includegraphics[height=1.40in,width=1.20in]{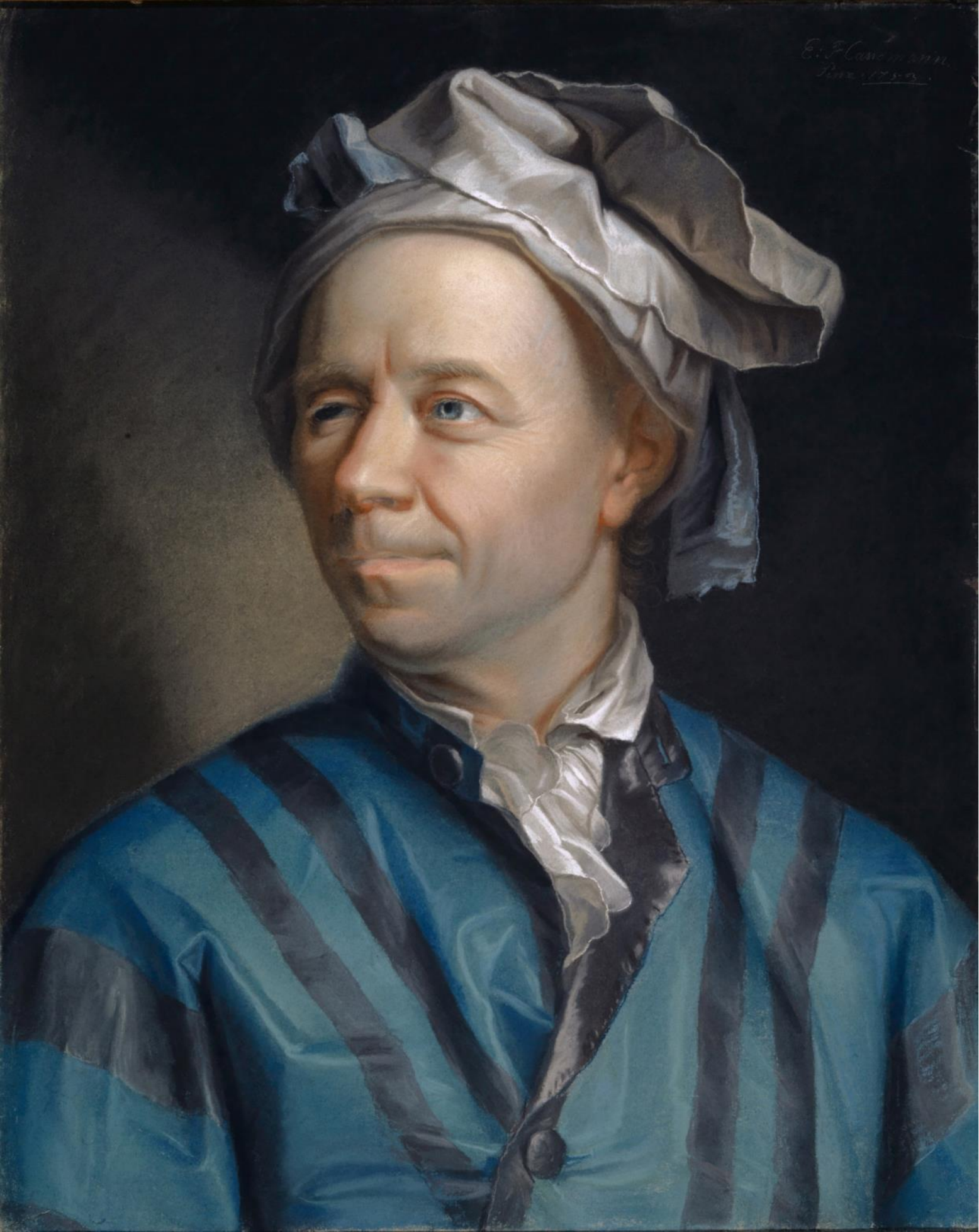}
\caption{\scriptsize \,\,
Leonhard Euler (15 April 1707 -- 18 September 1783)
formulated the Euler equations for fluid mechanics
which are among the first PDEs to be written down}
\label{figure-6a}
\end{minipage}
\begin{minipage}{0.485\textwidth}
\centering
\includegraphics[height=1.40in,width=1.1in]{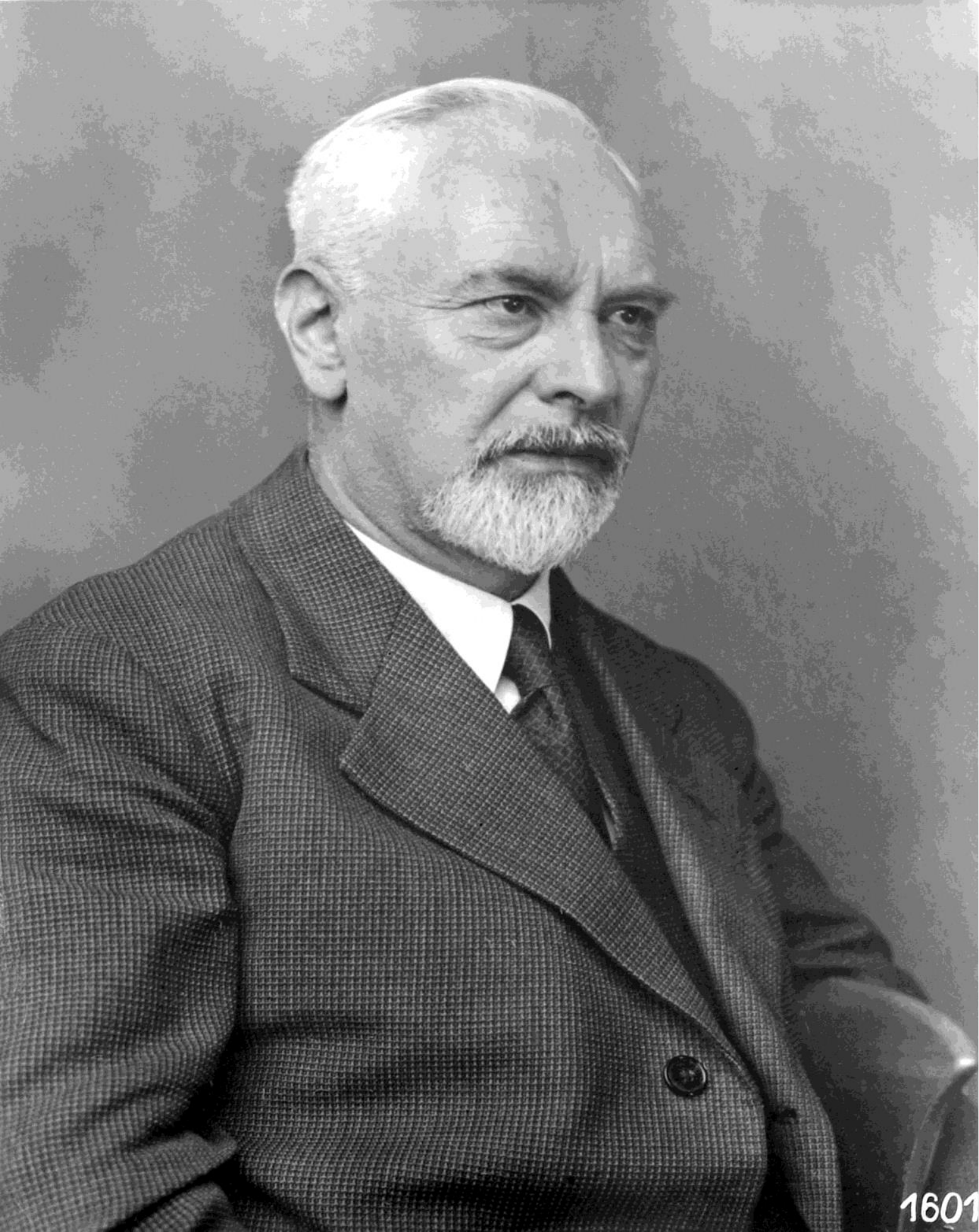}
\caption{\scriptsize Ludwig Prandtl (4 February 1875 -- 15 August 1953)
identified two oblique shock configurations when a steady uniform
supersonic gas flow passes a solid wedge via the shock polar analysis in 1936}
\label{figure-6b}
\end{minipage}
\end{figure}

\begin{figure}\label{figure-4}
\includegraphics[width=11.10cm,height=5.4cm]{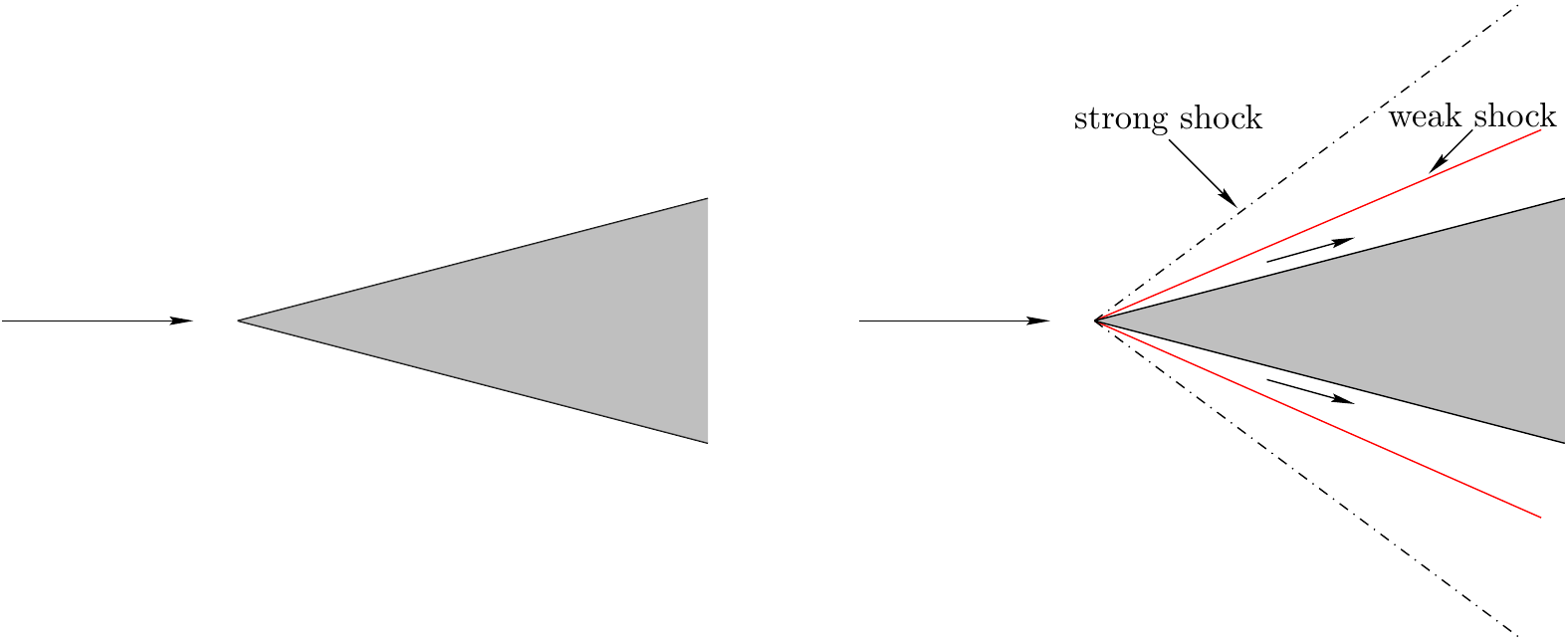}
\caption{\scriptsize Two steady solutions with shocks around the solid wedge
with angle $\theta_{\rm w}\in (0, \theta_{\rm w}^{\rm s})$ or even
$\theta_{\rm w} \in [\theta_{\rm w}^{\rm s}, \theta_{\rm w}^{\rm d})$}
\end{figure}

\begin{figure}
 \centering
\includegraphics[height=34mm]{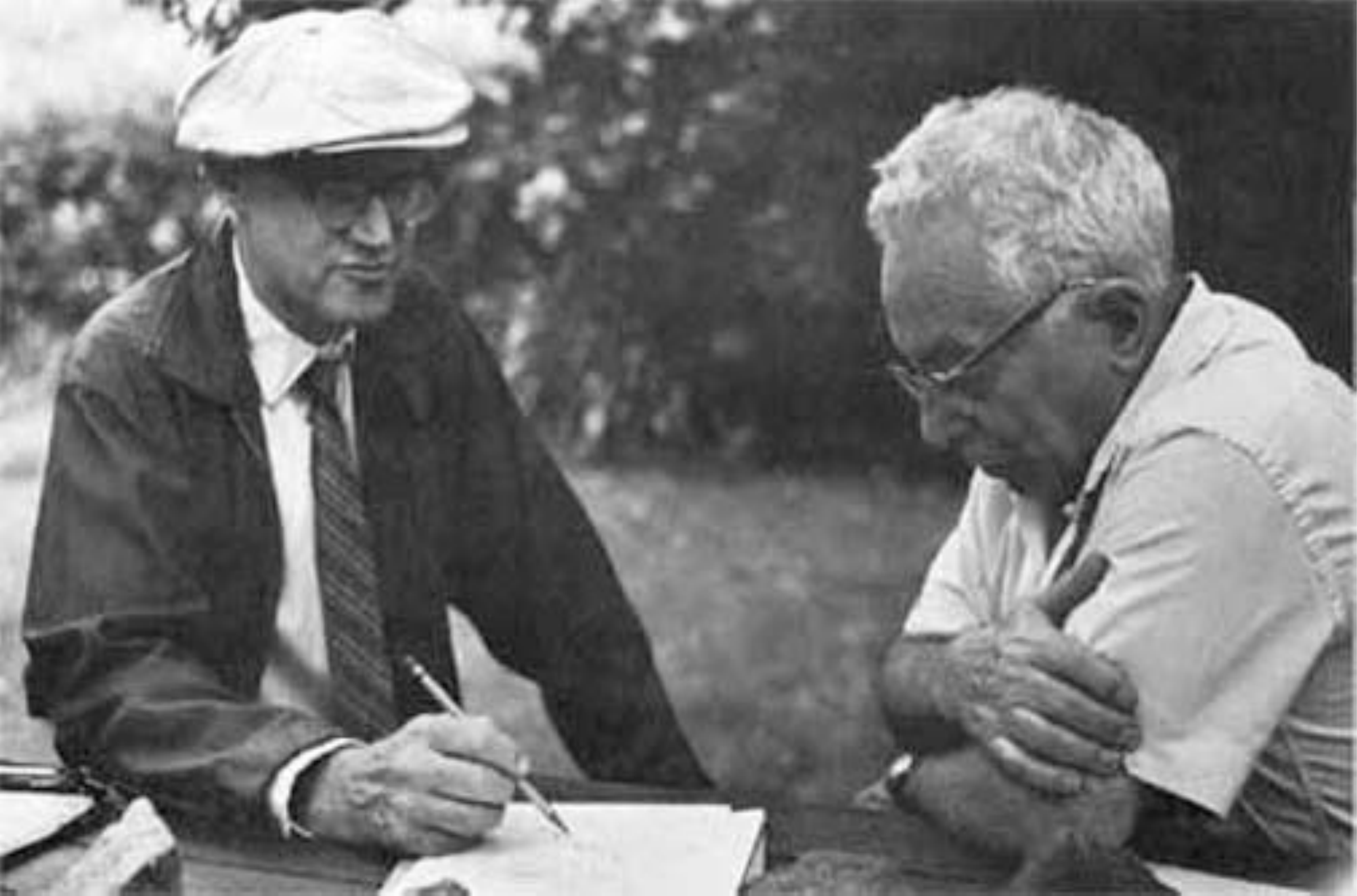}
\caption[]{\scriptsize Richard Courant (8 January 1888 -- 27 January 1972)
and Kurt Otto Friedrichs (28 September 1901 -- 31 December 1982);
their monumental book \cite{CFr}
has had great impact on
the development of the M-D theory of shock waves and nonlinear PDEs of hyperbolic/mixed types}
\label{Figure-5a}
\end{figure}

\medskip
Many fundamental transonic flow problems in fluid mechanics involve these nonlinear PDEs of mixed type.
One of them is a classical shock problem, in which an upstream steady uniform supersonic gas flow passes
a symmetric straight-sided solid wedge
\begin{equation}\label{wedge-1}
W:=\{\xx=(x_1,x_2)\in \R^2\,:\, |x_2|<x_1\tan \theta_{\rm w}, x_1>0\},
\end{equation}
whose (half-wedge) angle $\theta_{\rm w}$ is less than
the detachment angle $\theta^{\rm d}_{\rm w}$ ({\it cf}. Fig. 7).

Since this problem involves shocks, its global solution should be a weak solution of
equation \eqref{steady-potential-eq} or system  \eqref{steady-euler-syst}
 in the distributional sense (which admit shocks)\footnote{
Lax, P.~D.: Hyperbolic Systems of Conservation Laws and the Mathematical Theory of Shock Waves. CBMS-RCSAM,
No. 11. SIAM: Philadelphia, Pennsylvania, 1973.}
in the domain under consideration (see \cite{CF22}).
For example, for equation \eqref{steady-potential-eq}, a shock is a curve across which $\nabla\vphi$ is discontinuous.
If $\Lambda^+$ and $\Lambda^-(:=\Lambda\setminus \ol{\Lambda^+})$
are two nonempty open subsets of a domain $\Lambda\subset \R^2$, and $\CS:=\partial\Lambda^+\cap \Lambda$ is a $C^1$-curve
across which $\nabla\vphi$ has a jump, then $\vphi\in
C^1(\Lambda^{\pm}\cup \CS)\cap C^2(\Lambda^{\pm})$
is a global weak solution of \eqref{steady-potential-eq} in $\Lambda$ if and only
if $\vphi$ is in $W^{1,\infty}_{\rm loc}(\Lambda)$\footnote{A $W^{k,p}$ function, for $1\le p\le \infty$ and $k\ge 1$ integer, is a real-value function such that
itself and its (weak) derivatives up to order $k$ are all $L^p$ functions.}
and satisfies equation \eqref{steady-potential-eq} and the Rankine-Hugoniot conditions on $\CS$:
\begin{align}
\vphi_{\Lambda^+\cap \CS}=\vphi_{\Lambda^-\cap \CS},\qquad
\rho_B(|\nabla\vphi|^2)\nabla\vphi\cdot\nnu|_{\Lambda^+\cap \CS}
=\rho_B(|\nabla\vphi|^2)\nabla\vphi\cdot\nnu|_{\Lambda^-\cap \CS},\label{a-i}
\end{align}
where $\nnu$ is the unit normal to $\CS$ in the flow direction, {\it i.e.}, $\nabla\vphi\cdot \nnu|_{\Lambda^\pm\cap \CS}>0$.
A piecewise smooth solution with discontinuities satisfying \eqref{a-i}
is called an {\it entropy solution}
of \eqref{steady-potential-eq} if it satisfies the entropy condition:
{\it the density $\rho$ increases in the flow direction of
$\nabla\vphi_{\Lambda^+\cap \CS}$ across any discontinuity}.
Then such a discontinuity is called a {\it shock} (also see \cite{CFr}).

For this problem, there are two
configurations:
the weak oblique shock reflection with supersonic/subsonic downstream
flow (determined by the sonic angle $\theta_{\rm w}^{\rm s}$)
and the strong oblique shock reflection with subsonic downstream
flow -- both satisfy the entropy condition,
as discovered
by Prandtl ({\it cf}. Fig. \ref{figure-6b}).
The weak oblique shock is transonic with subsonic downstream flow for
$\theta_{\rm w}\in (\theta_{\rm w}^{\rm s}, \theta_{\rm w}^{\rm d})$,
while the weak oblique shock is supersonic with supersonic downstream flow for
$\theta_{\rm w}\in (0, \theta_{\rm w}^{\rm s})$.
However, the strong oblique shock is always transonic with subsonic downstream flow.
The question of physical admissibility of one
or both of the
strong/weak
shock reflection configurations
had been debated over the past eight decades since Courant-Friedrichs \cite{CFr}
and von Neumann \cite{vonN},
and has been better understood only recently
({\it cf.} \cite{CF22} and the references cited therein).
Two natural approaches for understanding this phenomenon are to examine whether these configurations
are stable under steady perturbations and to determine whether these configurations
are attainable
as large-time asymptotic states ({\it i.e.}, the {\it Prandtl-Meyer problem});
both approaches involve the analysis of nonlinear PDEs
\eqref{steady-potential-eq} or \eqref{steady-euler-syst} of mixed type.

Mathematically, the steady stability problem can be formulated
as a free boundary problem with the perturbed shock-front:
\begin{equation}\label{shock:1}
\sS=\{\xx \,: \, x_2 = \s(x_1),  \,x_1 \ge 0\} \qquad\quad\mbox{with $\s(0)=0\, $ and $\s(x_1)>0$ for $x_1>0$},
\end{equation}
as a free boundary (with the Rankine-Hugoniot conditions, say \eqref{a-i},
as free boundary conditions)
to determine the domain behind $\sS$:
\begin{equation}\label{domain:1}
\Omega=\{\xx\in \R^2\,:\,  b(x_1)<x_2<\s(x_1), x_1>0\}.
\end{equation}
and the downstream flow  in $\Omega$ for the nonlinear equation
\eqref{steady-potential-eq} or system \eqref{steady-euler-syst} of mixed elliptic-hyperbolic type,
where $x_2=b(x_1)$ is the perturbation of the flat wedge boundary $x_2=x_1\tan\theta_{\rm w}$.
Such a global solution of the free boundary problem
provides not only the global structural stability of
the steady oblique shock, but also a more detailed structure of the solution.

\begin{figure}
 \centering
\includegraphics[height=36mm]{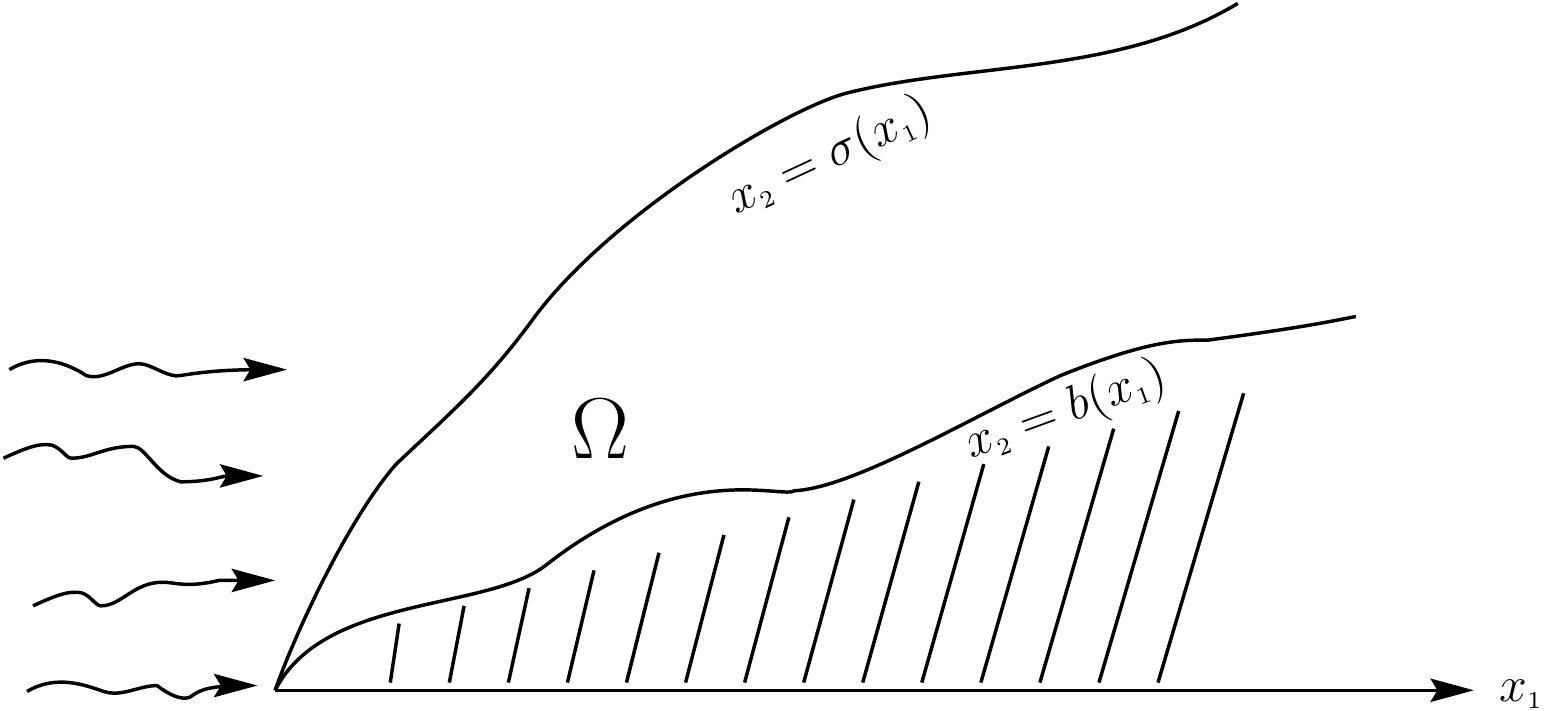}
\caption[]{\scriptsize The leading steady shock $x_2=\sigma(x_1)$ as a free boundary under the perturbation}
\label{Figure2a}
\end{figure}

Supersonic ({\it i.e.}, supersonic-supersonic) shocks
correspond to the case when $\theta_{\rm w}\in (0, \theta_{\rm w}^{\rm s})$,
which are shocks of weak strength.
The local stability of such shocks was first established in the 1960s.
The global stability and uniqueness of the supersonic oblique shocks
for both equation \eqref{steady-potential-eq} and system \eqref{steady-euler-syst}
have been solved for more general perturbations of both the upstream steady flow and
the wedge boundary,
even in $BV$\footnote{A $BV$ function is a real-valued functions whose total variation is bounded.},
by purely hyperbolic methods and techniques
({\it cf.} \cite{CF22} and the references cited therein).

\smallskip
For transonic ({\it i.e.}, supersonic-subsonic) shocks,
it has been proved that the oblique shock of weak strength is always stable under general steady perturbations.
However, the oblique shock of strong strength is stable only conditionally
for a certain class of steady perturbations that require
the exact match
of the steady perturbations near the wedge-vertex and the downstream condition at infinity, which reveals
one of the reasons why the strong oblique shock solutions have not been observed by experiment.
In these stability problems for transonic shocks,
the PDEs (or parts of the systems) are expected to be elliptic for global solutions in the domains determined
by the corresponding free boundary problems. That is, we solve an expected elliptic free boundary problem.
However, the earlier methods and approaches of elliptic free boundary problems do not directly apply to
these problems, such as
the variational methods, the Harnack inequality approach, and other elliptic methods/approaches.
The main reason is that the type of equations needs to be controlled before we can apply these methods, which
requires some strong {\it a priori} estimates. To overcome these difficulties, the global structure
of the problems is exploited, which allows us to derive certain properties of the solution
so that the type of equations and the geometry of the problem can be controlled. With this,
the free boundary problem as described above
has been solved by an iteration procedure.
See Chen-Feldman \cite{CF22} and the references cited therein for more details.

\begin{figure}
\begin{minipage}{0.49\textwidth}
\centering
\includegraphics[height=1.10in,width=1.90in]{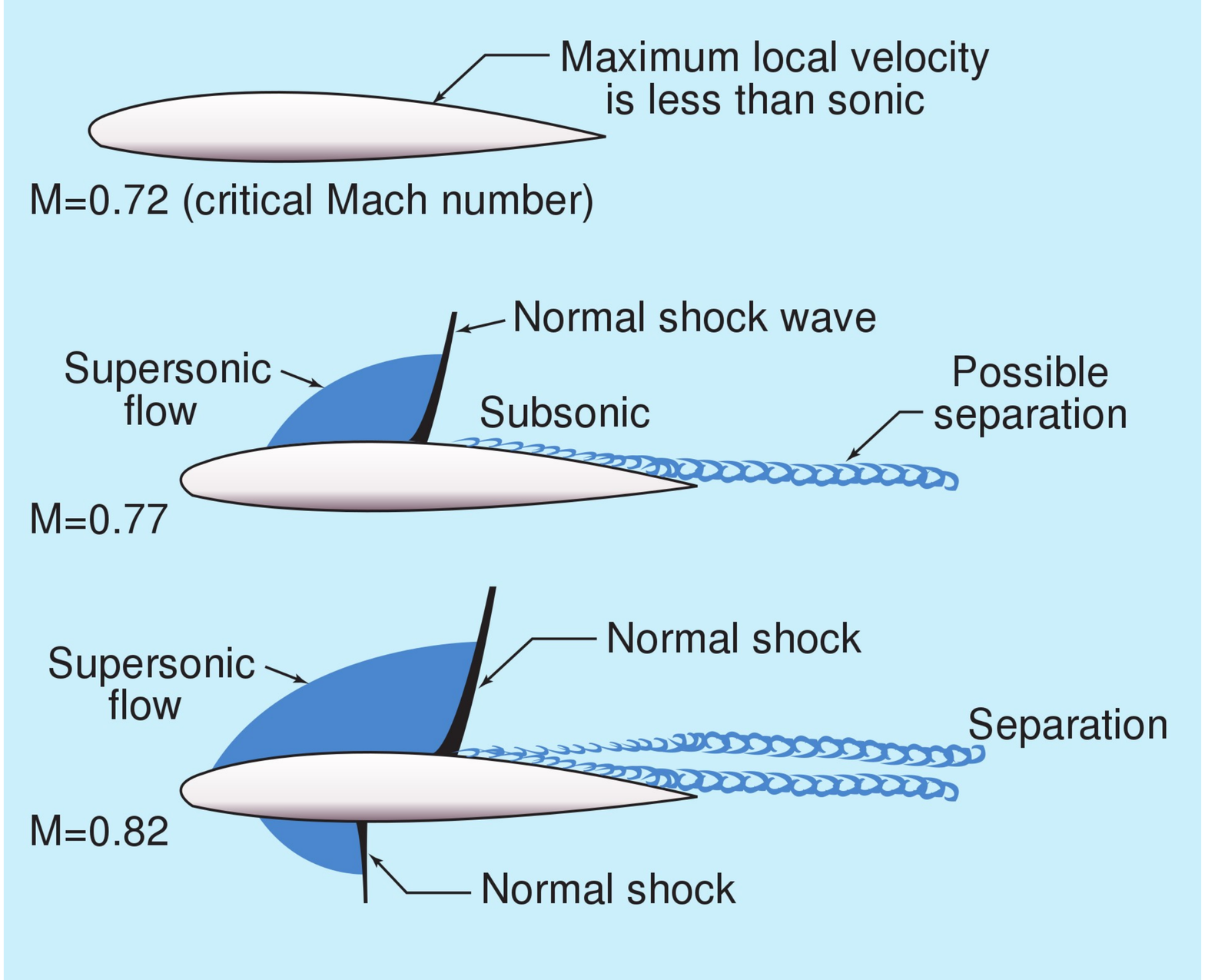}
\caption{\scriptsize Transonic flow patterns on an airfoil showing flow patterns at and above the critical Mach number}\label{fig-7}
\end{minipage}
\begin{minipage}{0.49\textwidth}
\centering
\includegraphics[height=1.10in,width=1.6in]{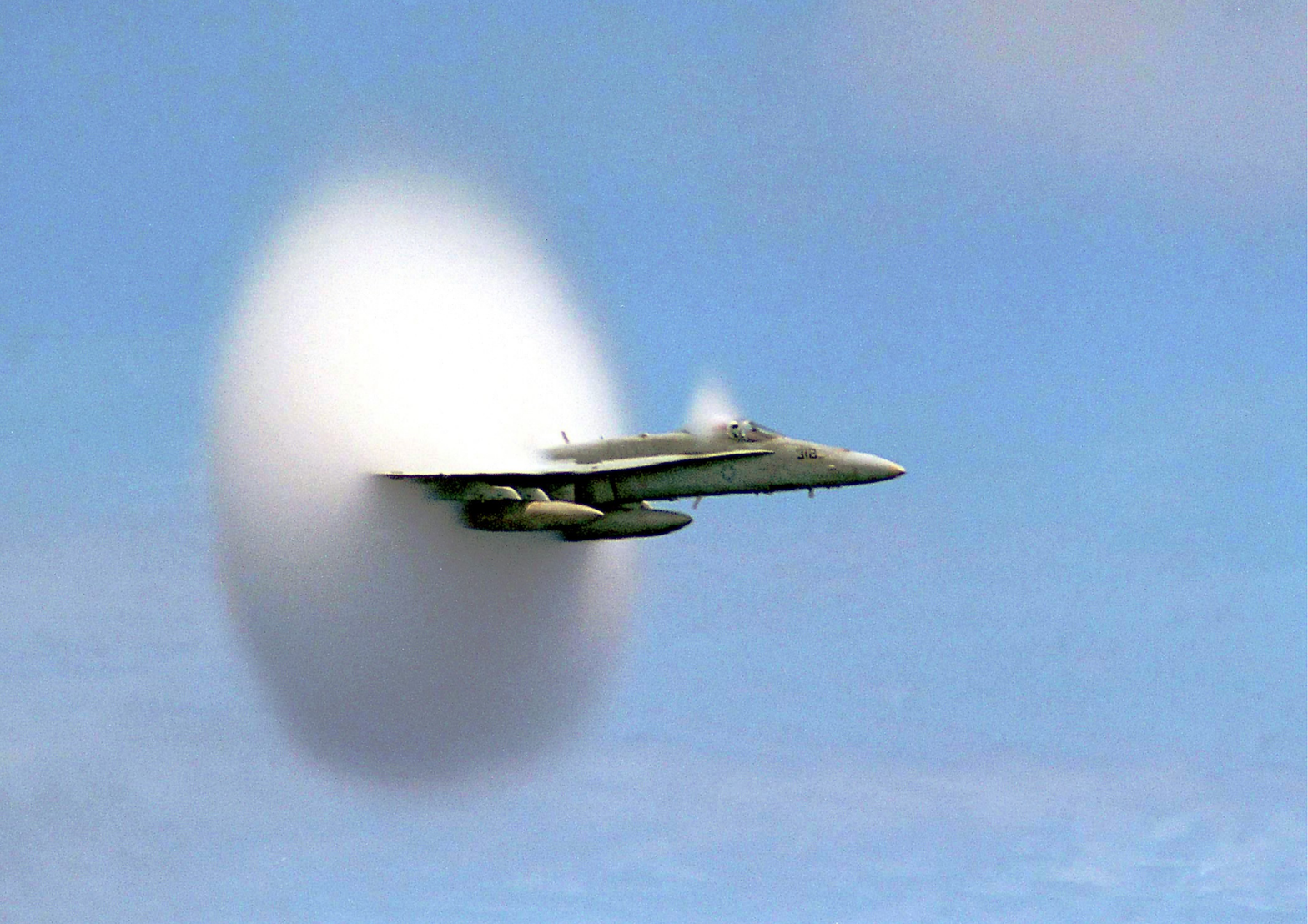}
\caption{
\scriptsize Aerodynamic condensation evidences of supersonic expansion fans around a transonic aircraft}\label{fig-8}
\end{minipage}
\end{figure}

\vspace{3pt}
When a subsonic flow passes through a {\it de Laval nozzle}, the flow may
form a supersonic bubble with a transonic shock
(see Fig. \ref{figure-3b});
the full understanding of how the geometry of the nozzle helps to create/stabilize/destabilize the transonic shock requires
a deep understanding of the nonlinear PDEs of mixed type.
Likewise, for the {\it Morawetz problem} for a steady subsonic
flow past an airfoil,
experimental results show
that a supersonic
bubble may be formed around the airfoil (see Figs. \ref{fig-7}--\ref{fig-8}), and  the flow behavior is
determined by the solution of the nonlinear PDE of mixed type.

Some fundamental
problems for transonic flow
posed in the 1950s--60s ({\it e.g.}, \cite{Bers,CF18,CFr,GM,vN1})
remain unsolved,
though some progress has been made in recent
years ({\it e.g.},
\cite{CF18,CF22,Da} and the references cited therein).

\section{Nonlinear PDEs of Mixed Type and Shock Reflection/Diffraction Problems in Fluid Mechanics and Related Areas}

In general, fluid flows are time-dependent.
We now describe how some longstanding M-D time-dependent fundamental
shock problems in fluid mechanics can naturally be formulated as problems
for nonlinear PDEs of mixed type,
through a prototype -- the {\it shock reflection-diffraction problem}.

When a planar shock
separating two constant states (0) and (1),
with constant velocities
and densities $\rho_0<\rho_1$
(state (0) is ahead or to the right of the shock, and state
(1) is behind the shock), moves in the flow direction ({\it i.e.}, $v_1>0$)
and hits a symmetric wedge \eqref{wedge-1} with (half-wedge) angle $\theta_{\rm w}$
head-on at time $t = 0$,
a reflection-diffraction process takes place for $t > 0$.
A fundamental question is which types of wave patterns of
shock reflection-diffraction configurations may be formed around the wedge.
The complexity of these
configurations was first reported
by Ernst Mach ({\it cf}. Fig. 12),
who observed two patterns of
shock reflection-diffraction configurations:
Regular reflection (two-shock configuration)
and Mach reflection (three-shock/one-vortex-sheet configuration),
as shown in Fig. 14 below\footnote{ M. Van Dyke:
{\it An Album of Fluid Motion}, The Parabolic Press: Stanford, 1982.}.
The issue remained dormant until the 1940s, when John von Neumann \cite{vN1,vonN} (also {\it cf}. Fig. 13),
as well as other mathematical/experimental
scientists ({\it cf.} \cite{BD,CF18,CFr,GM}
and the references cited therein),
began extensive research into all aspects of shock reflection-diffraction phenomena.
It has been found that the situations are much more complicated than
what Mach originally observed:
The shock reflection can be further
divided into more specific sub-patterns, and various other patterns of
shock reflection-diffraction configurations may occur, such as
the supersonic regular reflection, the subsonic regular reflection,
the attached regular reflection, the double Mach reflection,
the von Neumann reflection, and the Guderley reflection;
see \cite{BD,CF18,CFr,GM}
and the references cited therein (also see Figs. 14--\ref{fig:SubsonicRegularReflection-a} below).
Then the fundamental scientific issues include:
\begin{itemize}
\item[(i)] Structures of the shock reflection-diffraction configurations;

\item[(ii)] Transition criteria between the different patterns of the
configurations;

\item[(iii)] Dependence of the patterns upon the physical parameters such as the
wedge angle $\theta_{\rm w}$, the incident-shock-wave Mach number ({\it i.e.}, the strength of the incident shock),
and the
adiabatic exponent $\gamma>1$.
\end{itemize}

\begin{figure}
\begin{minipage}{0.502\textwidth}
\vspace{-10pt}
\centering
\includegraphics[height=1.70in,width=1.55in]{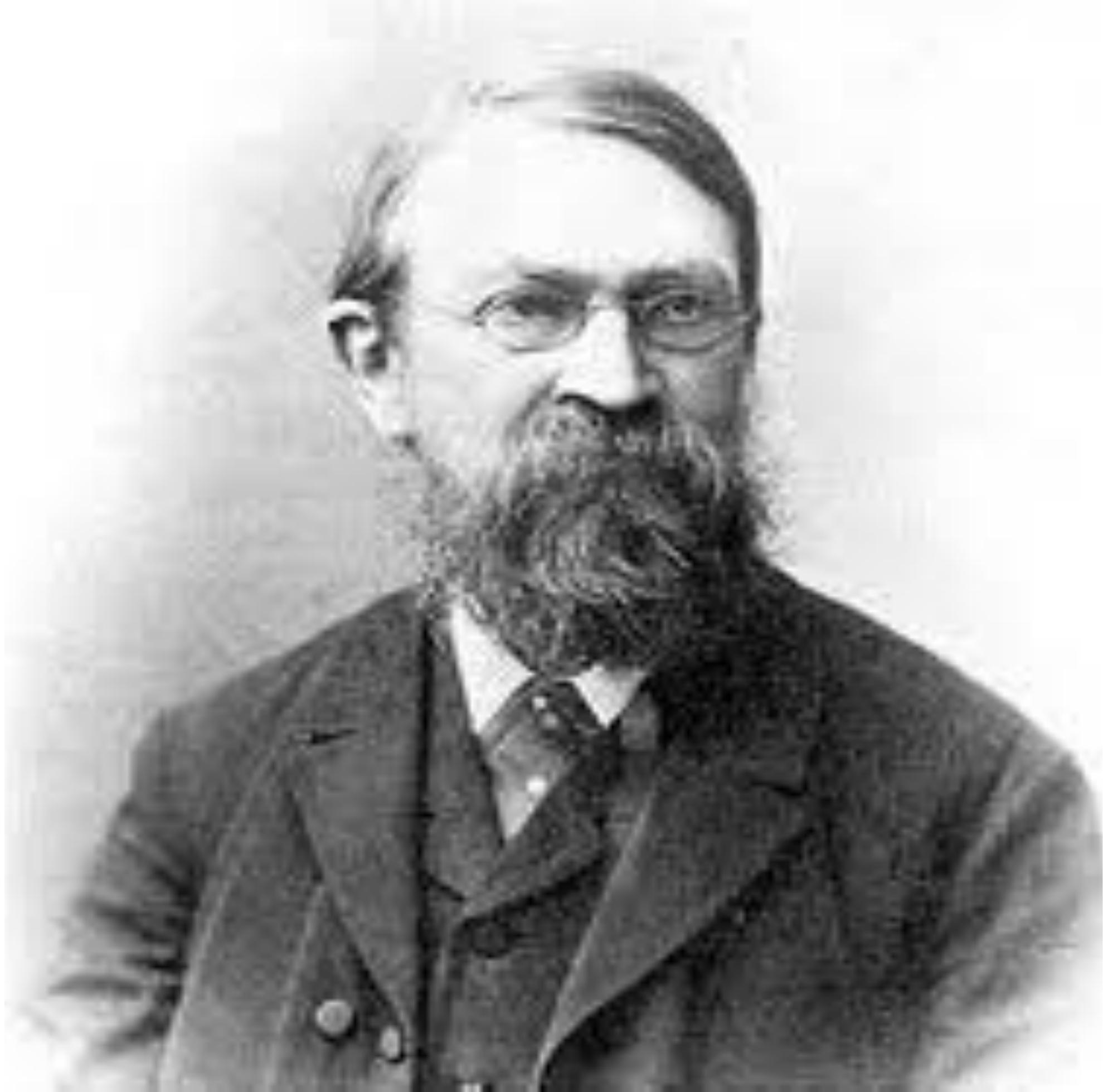}
\caption{\scriptsize \,\,
Ernst Waldfried Josef Wenzel Mach (18 February 1838 -- 19 February 1916)
who first observed the complexity of shock reflection-diffraction configurations (1878)}
\label{Mach}
\end{minipage}
\begin{minipage}{0.480\textwidth}
\centering
\includegraphics[height=1.60in,width=1.15in]{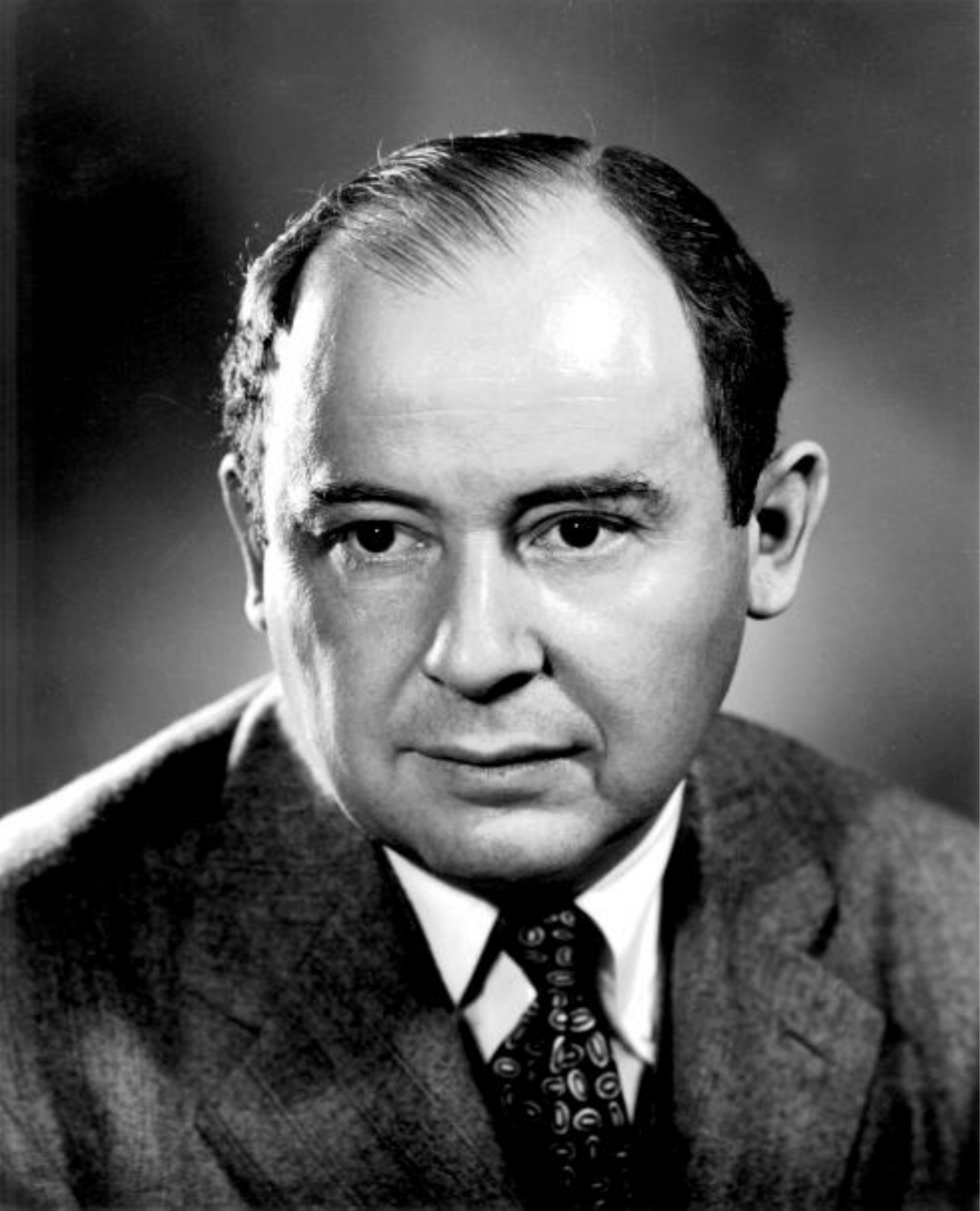}
\caption{\scriptsize John von Neumann (28 December 1903 -- 8 February 1957)
who proposed the sonic conjecture and the detachment conjecture for
shock reflection-diffraction configurations}
\label{Neumann}
\end{minipage}
\end{figure}

\noindent
In particular, several transition criteria between the different
patterns of shock reflection-diffraction configurations have been proposed,
including the {\it sonic conjecture} and the {\it detachment conjecture} by von Neumann
\cite{vN1} (also see \cite{BD,CF18}).

 \begin{figure}[h]\label{cog}
  \centering
  \includegraphics[width=0.18\textwidth]{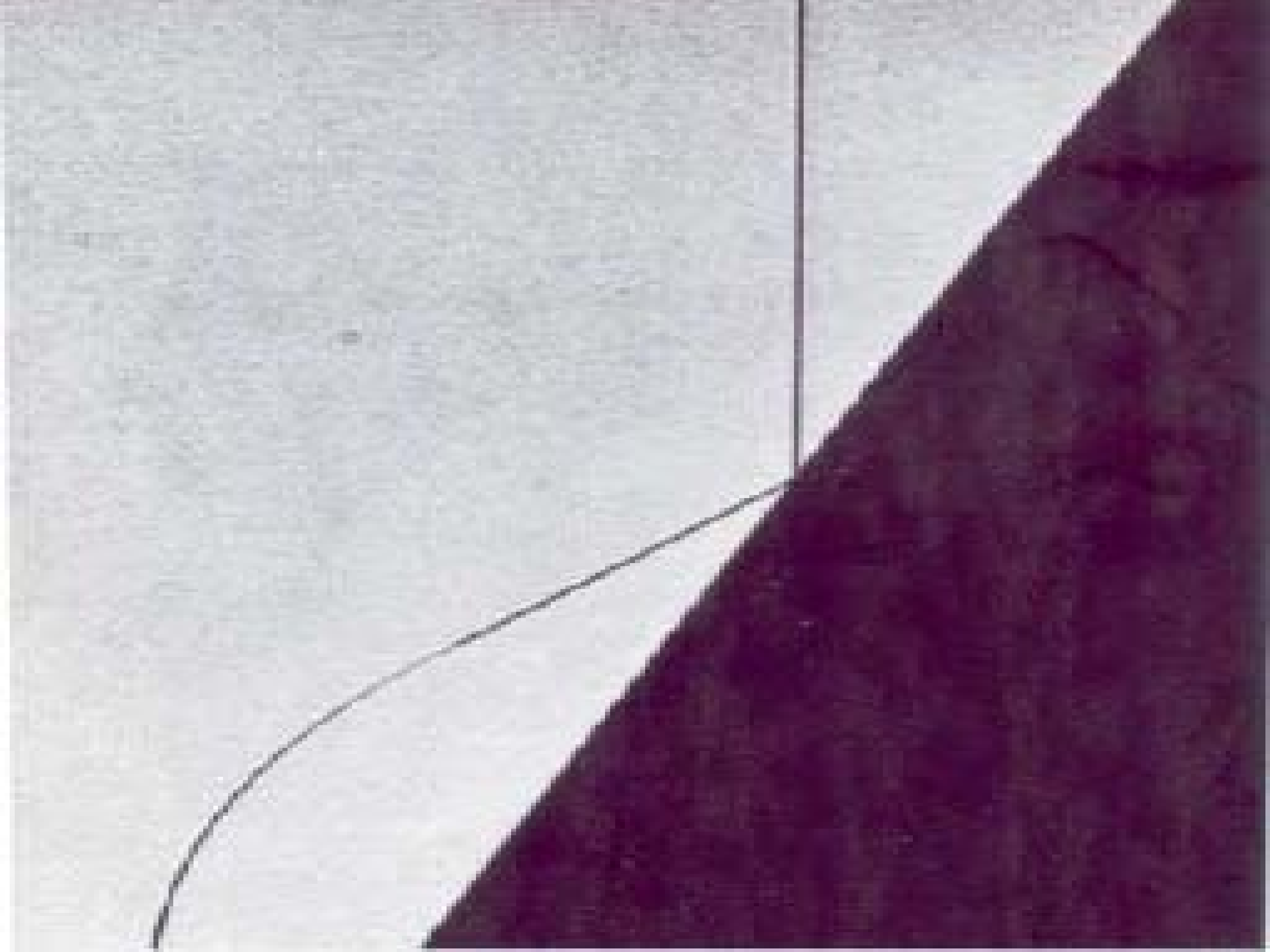} \qquad
  \includegraphics[width=0.18\textwidth]{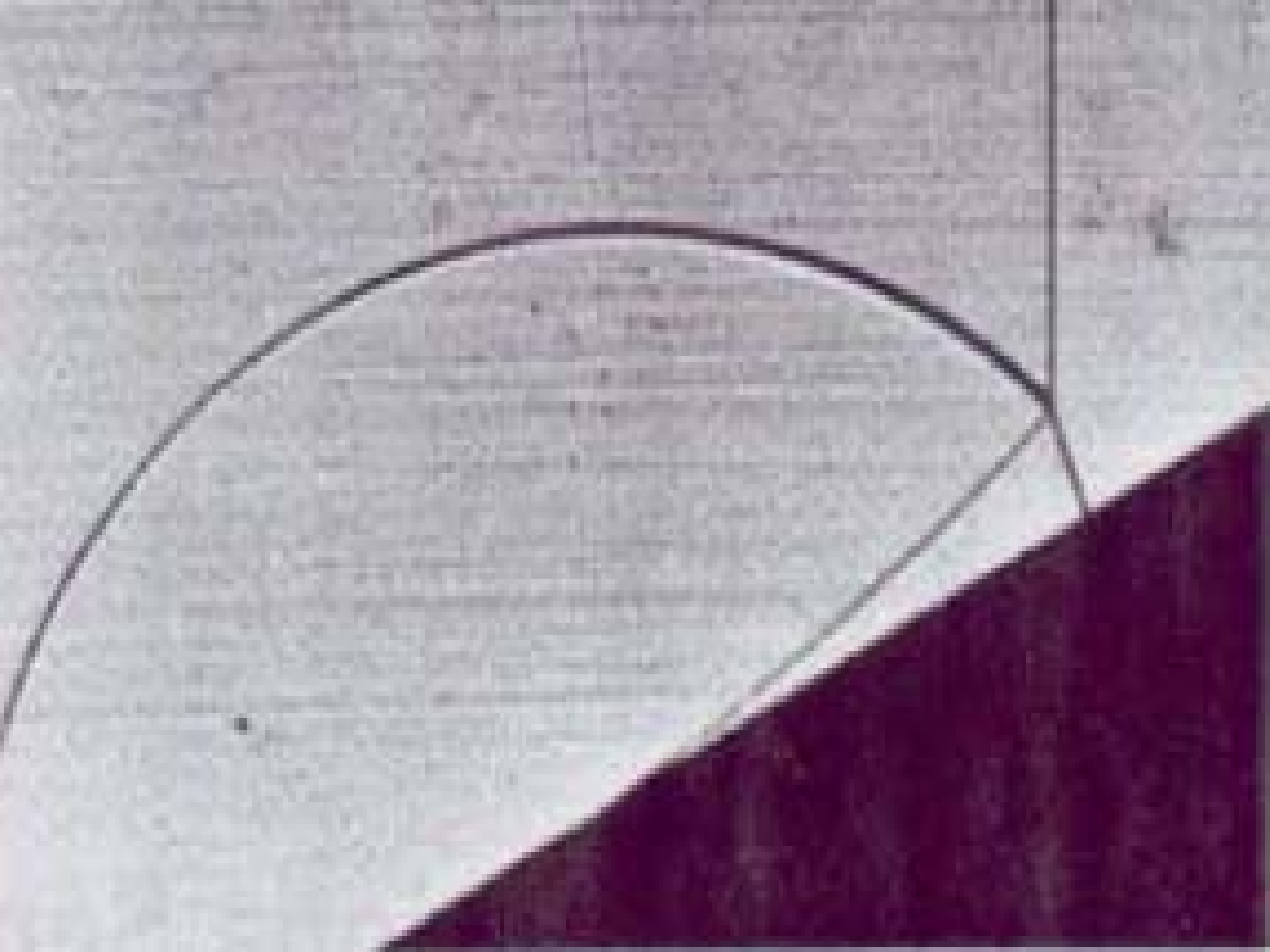}\qquad
  \includegraphics[width=0.18\textwidth]{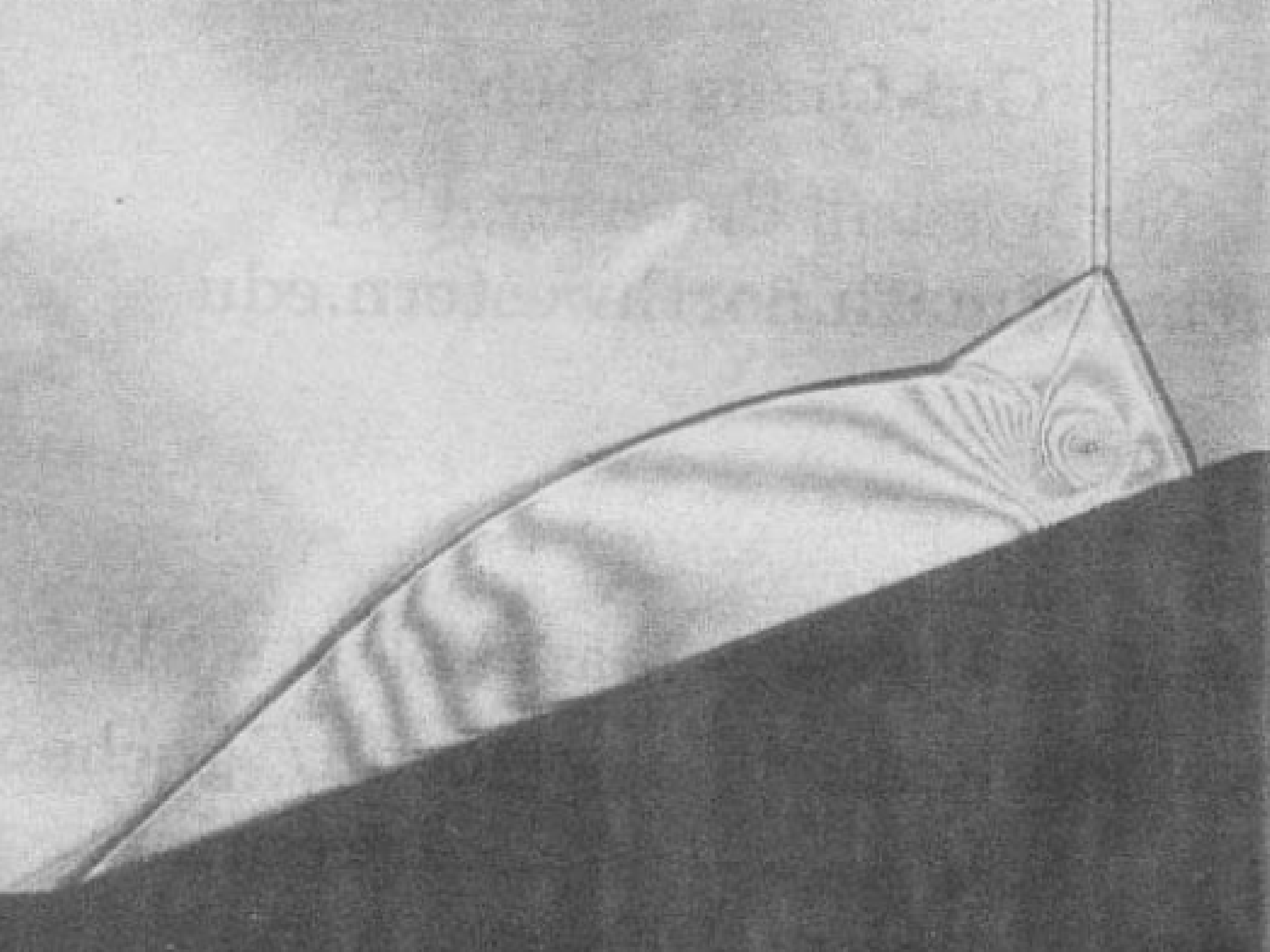}
  \caption[]{\scriptsize Three patterns of shock reflection-diffraction configurations}
\end{figure}

To present this more clearly, we focus now on
the Euler equations for
time-dependent compressible
potential flow, which consist of the conservation law
of mass and Bernoulli's law:
\begin{equation}\label{potential-t-eq}
\partial_t\rho +\divg(\rho\, \nabla\Phi)=0, \quad
\partial_t \Phi+{1\over 2}|\nabla\Phi|^2+{1\over \gamma-1}\rho^{\gamma-1}
=\frac{\rho_0^{\gamma-1}}{\gamma-1}
\qquad\quad\,\, \mbox{for $(t, {\bf x})\in \bR_+\times \bR^2$},
\end{equation}
after scaling, where
$\Phi$ is the time-dependent velocity potential ({\it i.e.}, $\vv=\nabla\Phi$ is the velocity).
Equivalently, system \eqref{potential-t-eq} can be reduced to the
nonlinear wave equation of second-order:
\begin{equation}\label{NonlinearWave-eq}
\del_t {\rho}(\del_t\Phi, \nabla_\x\Phi)+
\nabla_\x\cdot\big({\rho}(\del_t\Phi, \nabla_\x\Phi)\nabla_\x\Phi\big)=0,
\end{equation}
with
$
\rho(\del_t\Phi,\nabla_\x\Phi)=
\big(\rho_0^{\gamma-1}-(\gamma-1)(\del_t\Phi
+\frac{1}{2}|\nabla_\x\Phi|^2)\big)^{\frac{1}{\gamma-1}},
$
which is one of the original motivations for the extensive study of {\it nonlinear wave equations}.

\smallskip
Mathematically, the shock reflection-diffraction problem
is a $2$-D lateral Riemann problem for \eqref{potential-t-eq} or \eqref{NonlinearWave-eq}
in domain
$\R^2\setminus \overline{W}$ with
$\rho_0, \rho_1, v_1>0$ satisfying
\begin{equation}\label{ID-1}
\rho_1>\rho_0, \quad v_1=(\rho_1-\rho_0)\sqrt{\frac{2(\rho_1^{\gamma-1}-\rho_0^{\gamma-1})}{\rho_1^2-\rho_0^2}}.
\end{equation}

\begin{problem}[Shock Reflection-Diffraction Problem]\label{ibvp-c}
{\it
Piecewise constant initial data, consisting of state $(0)$ with velocity $\mathbf{v}_0=(0,0)$
and density $\rho_0>0$ on $\{x_1>0\}\setminus \overline{W}$
and state $(1)$ with velocity $\mathbf{v}_1=(v_1, 0)$ and density $\rho_1>0$
on $\{x_1 < 0\}$ connected by a shock at $x_1=0$,
are prescribed at $t = 0$, satisfying \eqref{ID-1}.
Seek a solution of the Euler system \eqref{potential-t-eq},
or Eq. \eqref{NonlinearWave-eq},
for $t\ge 0$ subject to these initial
data and the boundary condition $\nabla\Phi\cdot\nnu_{\rm w}=0$ on $\partial W$, where
$\nnu_{\rm w}$ is the unit outward normal to $\partial W$.
}
\end{problem}

\smallskip
{Problem \ref{ibvp-c}}
is invariant under scaling:
$(t, {\bf x}, \Phi)\rightarrow
(\tilde{t}, \tilde{\bf x}, \tilde{\Phi}(\tilde{t},\tilde{\bf x})):=(\alpha t, \alpha{\bf x}, \alpha\Phi(t,{\bf x}))$ for any $\alpha\neq 0$,
and thus it admits  self-similar solutions in the form of
\begin{equation}\label{4.6}
\Phi(t, {\bf x})=t\phi(\xxi)\qquad\quad \text{for}\;\;\xxi=\frac{{\bf x}}{t}.
\end{equation}
Then the pseudo-potential function
$\vphi(\xxi)=\phi(\xxi)-\frac 12|\xxi|^2$
satisfies
the following equation:
\begin{equation}\label{2-1}
{\rm div}\big(\rho_B(|\D\vphi|^2,\vphi)\D\vphi\big)+2\rho_B(|\D\vphi|^2,\vphi)=0
\end{equation}
with
$\rho_B(|\D\vphi|,\vphi)=
\bigl(\rho_0^{\gamma-1}- (\gamma-1)(\frac{1}{2}|\D\vphi|^2+\varphi)\bigr)^{\frac{1}{\gamma-1}},
$
where
the divergence ${\rm div}$ and gradient $\D$ are with respect to $\xxi\in \R^2$.
Define the pseudo-sonic speed $c=c(|\D\varphi|,\varphi)$ by
\begin{equation}\label{c-through-density-function}
c^2(|\D\varphi|,\varphi)=
\rho^{\gamma-1}(|\D\varphi|^2,\varphi)
=B_0-(\gamma-1)\big(\frac{1}{2}|\D\varphi|^2+\varphi\big).
\end{equation}
Eq. \eqref{2-1} is of mixed elliptic-hyperbolic type -- It is
\begin{itemize}
\item strictly {\it elliptic} if $|\D\varphi|<c(|\D\varphi|, \varphi)$
(pseudo-subsonic),

\smallskip
\item strictly {\it hyperbolic} if $|\D\varphi|>c(|\D\varphi|, \varphi)$
(pseudo-supersonic).
\end{itemize}
The transition boundary between
the pseudo-supersonic and pseudo-subsonic phases is $|\D\varphi|=c(|\D\varphi|, \varphi)$ ({\it i.e.},
$|\D\varphi|=\sqrt{\frac{2}{\gamma+1}\big(B_0-(\gamma-1)\varphi\big)}$\,),
a degenerate set of the solution of equation \eqref{2-1}, which is {\it a priori} unknown
and more delicate than that of equation \eqref{steady-potential-eq}.

One class of solutions of (\ref{2-1})
is that of {\em constant states} that are the solutions
with constant velocity $\mathbf{v}_*\in \mathbb{R}^2$.
Then
the pseudo-potential of a constant state satisfies
$\D\varphi=\mathbf{v}_*-\xxi$ so that
\begin{equation}\label{constantStatesForm}
\varphi(\xxi)=-\frac 12|\xxi|^2+\mathbf{v}_*\cdot\xxi +C,
\end{equation}
where $C$ is a constant. For this $\varphi$,
the density $\rho$ and sonic
speed $c=\rho^{(\gamma-1)/2}$ are positive constants, independent of $\xxi$.
Then, from
\eqref{constantStatesForm},
the ellipticity condition for the constant state is
$
|\xxi -\mathbf{v}_*|<c.
$
Thus, for a constant state $\mathbf{v}_*$,
equation (\ref{2-1})
is {\it elliptic} inside the {\em sonic circle},
with center $\mathbf{v}_*$ and radius $c$, and {\it hyperbolic} outside this circle.
Moreover, if the density $\rho$ is a constant, then the solution is a constant state;
that is, the corresponding
pseudo-potential $\vphi$ is of form \eqref{constantStatesForm}.

Problem \ref{ibvp-c} involves transonic shocks so that its global solution should be a weak solution of
equation \eqref{2-1}
in the distributional sense
within the domain in the $\xxi$--coordinates (see \cite{CF22}).
If $\Lambda^+$ and $\Lambda^-(:=\Lambda\setminus \ol{\Lambda^+})$
are two nonempty open subsets of a domain $\Lambda\subset \R^2$, and $\CS:=\partial\Lambda^+\cap \Lambda$ is
a $C^1$-curve with normal $\nnu$
across which $\D\vphi$ has a jump, then $\vphi\in
C^1(\Lambda^{\pm}\cup \CS)\cap C^2(\Lambda^{\pm})$
is a global entropy solution of \eqref{2-1} in $\Lambda$ with $\CS$ as a shock
{\it if and only if} $\vphi$ is in $W^{1,\infty}_{\rm loc}(\Lambda)$
and satisfies equation \eqref{2-1}, the Rankine-Hugoniot conditions on $\CS$:
\vspace{-4mm}
\begin{align}
&\vphi_{\Lambda^+\cap \CS}=\vphi_{\Lambda^-\cap \CS},\label{1-i}\\
&\rho(|\D\vphi|^2, \vphi)\D\vphi\cdot\nnu|_{\Lambda^+\cap \CS}
=\rho(|\D\vphi|^2, \vphi)\D\vphi\cdot\nnu|_{\Lambda^-\cap \CS},\label{1-h}
\end{align}
and the entropy condition: {\it the density $\rho$ increases in the pseudo-flow direction of
$\D\varphi_{\Lambda^+\cap \CS}$ across any discontinuity}.

We now show how such solutions of the nonlinear PDE \eqref{2-1} of mixed elliptic-hyperbolic type
in self-similar coordinates $\xxi
=\frac{{\bf x}}{t}$
can be constructed.

First, by the symmetry of the problem with respect to the $\xi_1$--axis,
it suffices for us to focus only on the upper half-plane
$\{\xi_2>0\}$ and prescribe the slip boundary condition: $\D\varphi \cdot \bn_{\rm sym}=0$
on the symmetry line $\Gamma_{\rm sym}:=\{\xi_2=0\}$ for the interior unit normal $\bn_{\rm sym}=(0,1)$.
Then {Problem \ref{ibvp-c}} can be reformulated as a boundary value problem
in the unbounded domain:
$$
\Lambda:=\R^2_+\setminus\{\xxi\,:\,|\xi_2|\le \xi_1 \tan\theta_{\rm w}, \xi_1>0\}
$$
in the self-similar coordinates $\xxi=(\xi_1,\xi_2)$, where $\mR^2_+:=\mR^2\cap\{\xi_2>0\}$.

\begin{problem}[Boundary Value Problem]\label{bvp-c}
{\it Seek a
solution $\varphi$ of equation \eqref{2-1}
in the self-similar
domain $\Lambda$ with the slip boundary condition{\rm :}
$\D\varphi\cdot\nnu|_{\partial\Lambda}=0$ for the interior unit normal $\nnu$
on $\partial\Lambda$,
and the asymptotic boundary condition at infinity{\rm :}
$$
\varphi\longrightarrow\bar{\varphi}=
\begin{cases} \varphi_0 \qquad\mbox{for}\,\,\,
                         \xi_1>\xi_1^0, \xi_2>\xi_1 \tan\theta_{\rm w},\\
              \varphi_1 \qquad \mbox{for}\,\,\,
                          \xi_1<\xi_1^0, \;\xi_2>0,
\end{cases}
\qquad \mbox{when $|\xxi|\to \infty$,}
$$
where $\pSi_0=-\frac{1}{2}|\xxi|^2$ and $\pSi_1=-\frac{1}{2}|\xxi|^2+v_1(\xi_1-\xi^0_1)$ with
$\xi_1^0
=\frac{\rho_1v_1}{\rho_1-\rho_0},
$
which is the location of the incident shock $\CS_0=\{\xi_1=\xi_1^0\}\cap\Lambda$
determined by the Rankine-Hugoniot conditions
\eqref{1-i}--\eqref{1-h}
between states $(0)$ and $(1)$ on $\CS_0$.
}
\end{problem}

\begin{figure}
 \centering
\includegraphics[height=37mm]{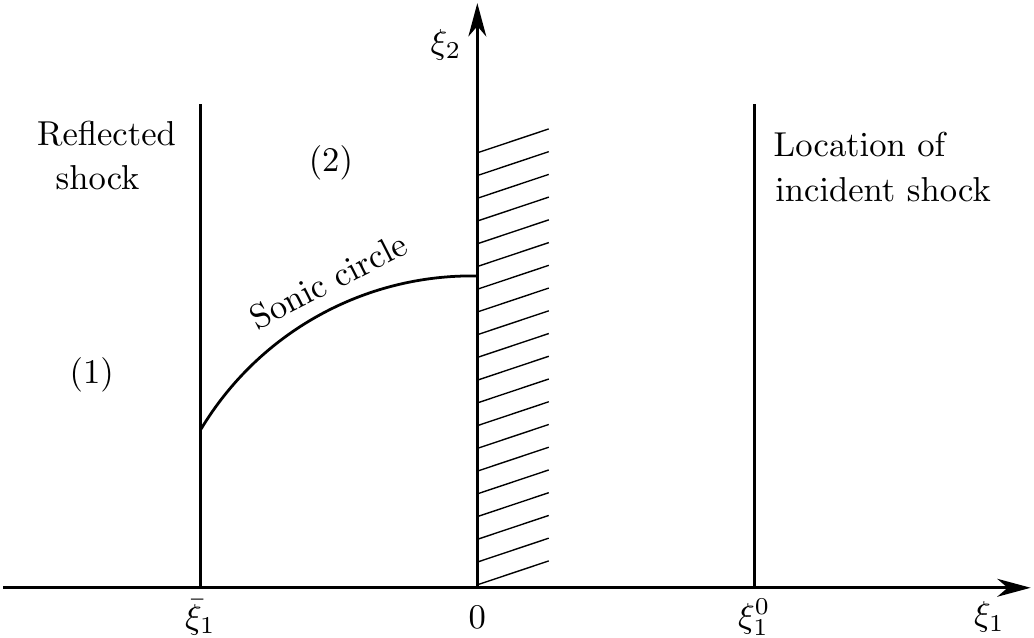}
\caption[]{\scriptsize Normal reflection configuration}
\label{NormReflFigure}
\end{figure}

\smallskip
The simplest case is when  $\theta_{\rm w}=\frac\pi 2$,
which is called {\em normal reflection}; see Fig. \ref{NormReflFigure}.
In this case, the incident shock normally reflects from the flat wall
to become
the flat reflected shock
$\xi_1=\bar{\xi}_1<0$.

When  $\theta_{\rm w}\in (0, \frac\pi 2)$,
a necessary condition for the existence
of a regular reflection solution whose configurations are as shown
in Figs. \ref{fig:RegularReflection}--\ref{fig:SubsonicRegularReflection-a}
is the existence of the uniform state (2) with pseudo-potential $\varphi_2$ at $P_0$, determined by
the following three conditions at $P_0$:
\begin{equation}\label{condState2}
\D\varphi_2\cdot\bn_{\rm w}=0, \,\,\,\,\varphi_2=\varphi_1,\,\,\,\,\r(|\D\varphi_2|^2,\varphi_2)\D\varphi_2\cdot\bn_{\CS_1}=
\rho_1\D\varphi_1\cdot\bn_{\CS_1} \qquad
\mbox{for $\bn_{\mathcal{S}_1}=\frac{\D(\varphi_1-\varphi_2)}{|\D(\varphi_1-\varphi_2)|}$}
\end{equation}
across the flat shock $\CS_1=\{\varphi_1=\varphi_2\}$ that separates state (2) from state (1) and satisfies
the entropy condition: $\rho_2>\rho_1$.
These conditions lead to the system of algebraic equations \eqref{condState2}
for the constant velocity $\mathbf{v}_2$
and
density $\rho_2$ of state (2).
For any fixed densities $0<\rho_0<\rho_1$ of states (0) and (1),
there exist a sonic angle $\theta_{\rm w}^{\rm s}$ and a detachment angle
$\theta_{\rm w}^{\rm d}$
satisfying
$$
0<\theta_{\rm w}^{\rm d}<\theta_{\rm w}^{\rm s}<\frac{\pi}{2}
$$
such that the algebraic system \eqref{condState2} has two solutions for
each $\theta_{\rm w}\in (\theta_{\rm w}^{\rm d}, \frac{\pi}{2})$, which become equal when
$\theta_{\rm w}=\theta_{\rm w}^{\rm d}$.
Thus, for each $\theta_{\rm w}\in (\theta_{\rm w}^{\rm d}, \frac{\pi}{2})$, there exist two
states (2), called weak and strong, with densities  $0<\rho_1<\rho_2^{\rm weak}<\rho_2^{\rm strong}$ (the entropy condition).
The weak state (2) is supersonic at the reflection point $\PtIncW$ for
$\theta_{\rm w}\in (\theta_{\rm w}^{\rm s}, \frac{\pi}{2})$,
sonic for $\theta_{\rm w}=\theta_{\rm w}^{\rm s}$,
and subsonic for $\theta_{\rm w}\in (\theta_{\rm w}^{\rm d}, \hat\theta_{\rm w}^{\rm s})$
for some $\hat\theta_{\rm w}^{\rm s}\in(\theta_{\rm w}^{\rm d}, \theta_{\rm w}^{\rm s}]$.
The strong state (2) is always subsonic at $\PtIncW$ for all
$\theta_{\rm w}\in (\theta_{\rm w}^{\rm d}, \frac{\pi}{2})$.

There had been a long debate on determining which of the two states (2) for
$\theta_{\rm w}\in (\theta_{\rm w}^{\rm d}, \frac{\pi}{2})$,
 weak or strong, is physical for the local theory; see
\cite{BD,CF22,CFr}.
Indeed, it has been shown in Chen-Feldman \cite{CF2010,CF22}
that the weak shock reflection-diffraction configuration tends to the unique normal
reflection in Fig. \ref{NormReflFigure}, but the strong one
does not,
when the wedge angle $\theta_{\rm w}$ tends to $\frac{\pi}{2}$.
The strength of the corresponding
reflected shock near $P_0$ in the
weak shock reflection-diffraction configuration
is relatively weak,
compared to the other shock given by the strong state (2).
From now on, for the given wedge angle
$\theta_{\rm w}\in (\theta_{\rm w}^{\rm d}, \frac{\pi}{2})$,
state (2) represents the unique weak state (2)
and $\varphi_2$ is its pseudo-potential.

\smallskip
If the weak state (2) is supersonic, the speeds of propagation of the solution
are finite, and state (2) is
determined completely by the local information: state (1),
state (0), and the location of point $P_0$. That is, any information
from the reflection-diffraction domain,
particularly the disturbance at corner $P_3$,
cannot travel towards the reflection point $P_0$.
However, if
it is subsonic, the information can reach $P_0$ and interact with
it, potentially
altering the subsonic reflection-diffraction configuration.
This argument motivated the following conjectures by
von Neumann in \cite{vN1} (also see \cite{BD,CF18}):

\begin{figure}
\centering
\begin{minipage}{0.470\textwidth}
\centering
\includegraphics[height=1.5in,width=2.1in]{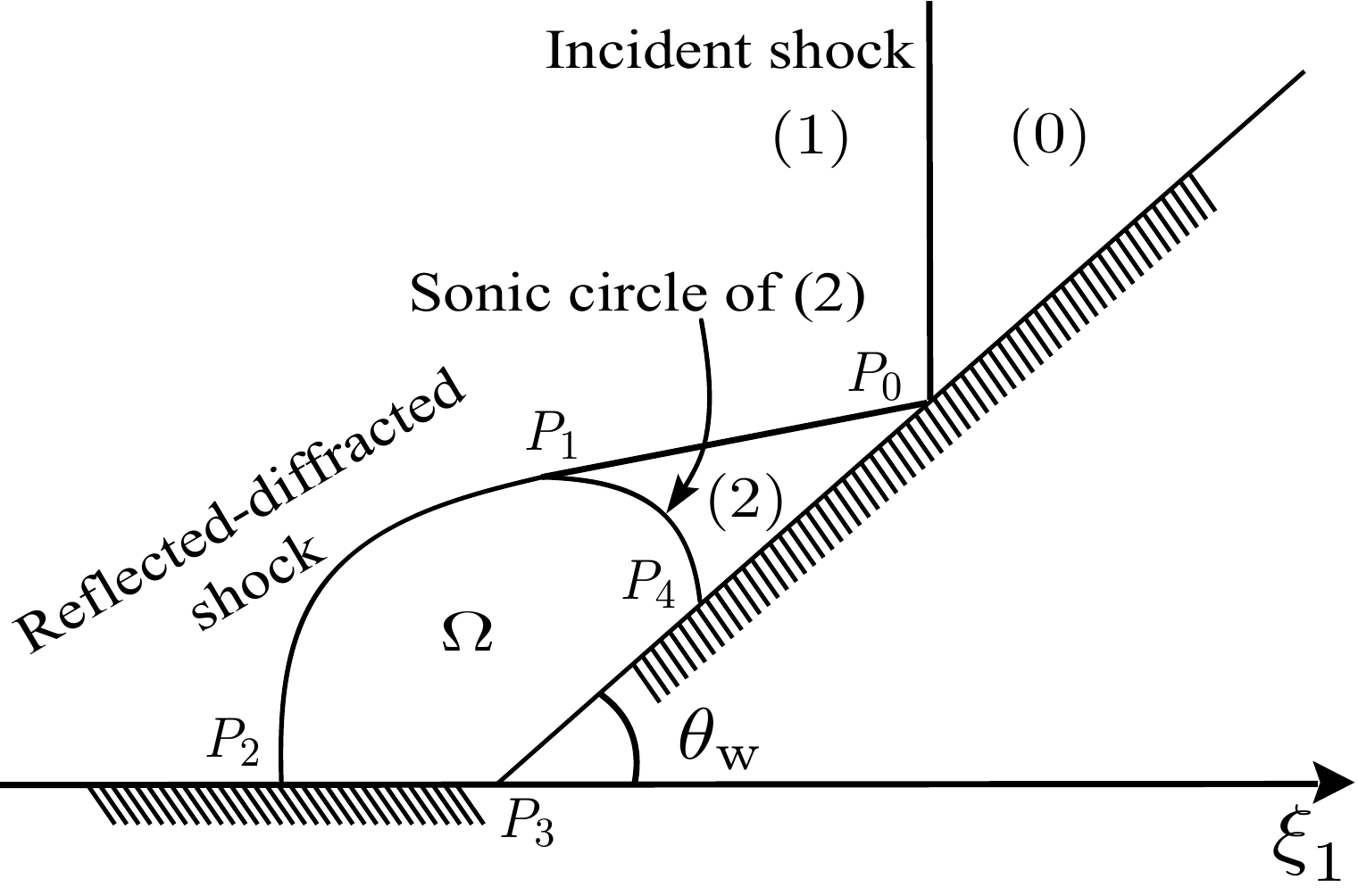}
\caption{\scriptsize \,\, Supersonic regular reflection-diffraction configuration \cite{CF18}}
\vspace{-2pt}
\label{fig:RegularReflection}
\end{minipage}
\,\,
\begin{minipage}{0.470\textwidth}
\centering
\includegraphics[height=1.5in,width=2.1in]{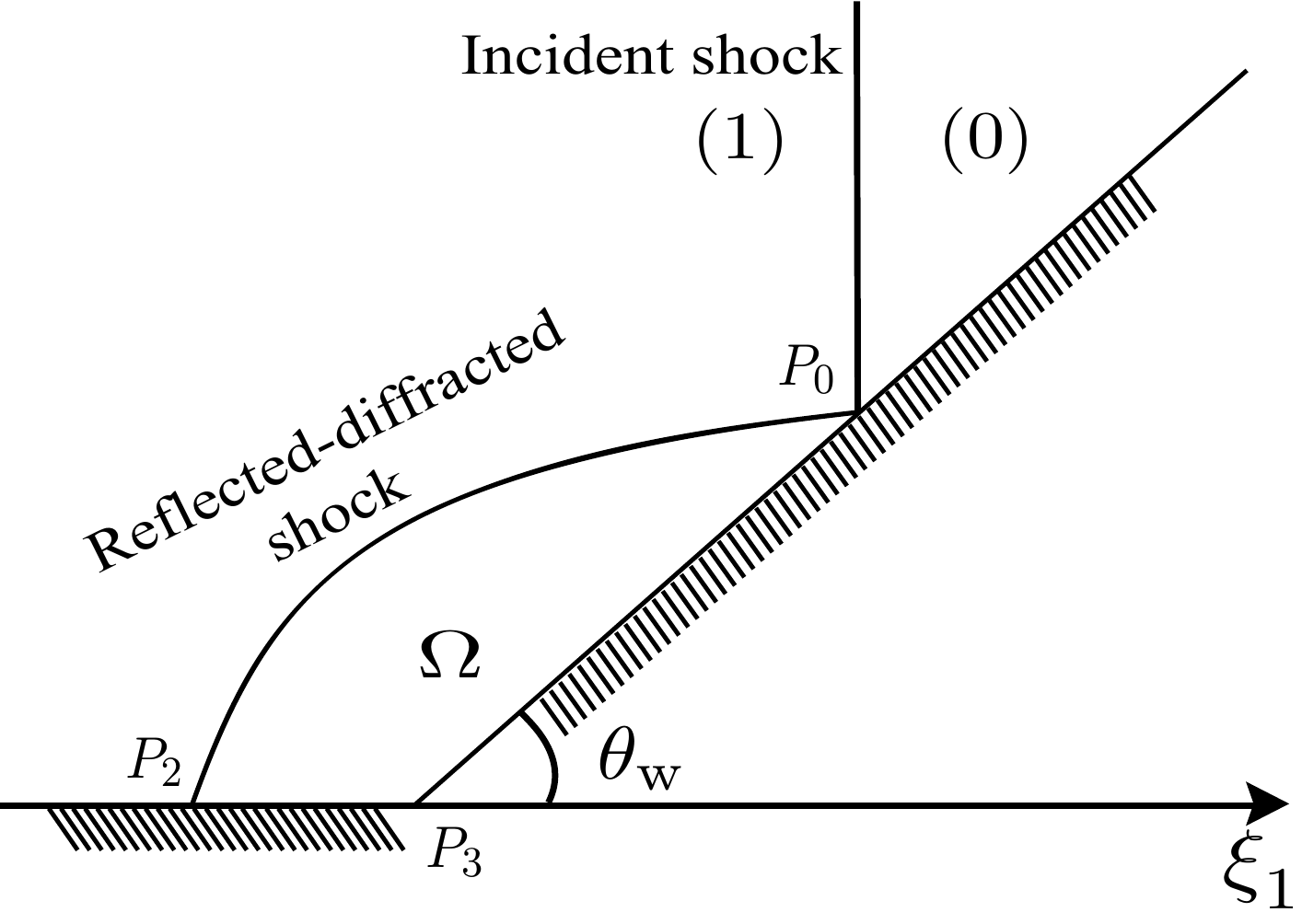}
\caption{\scriptsize \,\, Subsonic regular reflection-diffraction configuration \cite{CF18}}
\vspace{-2pt}
\label{fig:SubsonicRegularReflection}
\end{minipage}
\vspace{-7pt}
\end{figure}

\smallskip
{\bf The von Neumann Sonic Conjecture}:
{\em There exists a supersonic regular shock reflection-diffraction
configuration when
$\theta_{\rm w}\in (\theta_{\rm w}^{\rm s}, \frac{\pi}{2})$
for $\theta_{\rm w}^{\rm s}>\theta_{\rm w}^{\rm d}$.
That is,
the supersonicity of the weak state {\rm (2)} implies the existence
of a supersonic regular reflection
solution, as shown in {\rm Fig. \ref{fig:RegularReflection}.}}

\smallskip
Another conjecture is that the global regular shock reflection-diffraction
configuration is still possible whenever the local regular reflection at the reflection
point is possible:

\smallskip
{\bf The von Neumann Detachment Conjecture}:
{\em There exists a subsonic regular shock reflection-diffraction configuration for
any wedge angle $\theta_{\rm w}\in (\theta_{\rm w}^{\rm d}, \theta_{\rm w}^{\rm s})$.
That is, the existence of subsonic weak state {\rm (2)} beyond the sonic angle
implies the existence
of a subsonic regular reflection solution,
as shown in
{\rm Fig. \ref{fig:SubsonicRegularReflection}}}.

\smallskip
State (2)
determines the
straight shock $\CS_1$ and the sonic arc $\Gamma_{\rm sonic}:=P_1P_4$ when state (2) is supersonic at $P_0$,
and the slope of $\Shock$ at $P_0$ (arc $\Gamma_{\rm sonic}$ on the boundary of $\Omega$ becomes a corner point $P_0$)
when state (2) is subsonic at $P_0$.
Thus, the unknowns are the domain $\Omega$ (or equivalently, the curved part of the reflected-diffracted shock $\Shock$)
and the pseudo-potential $\varphi$ in $\Omega$.
Then, from \eqref{1-i}--\eqref{1-h}, in order to construct a solution of Problem \ref{bvp-c}
for the supersonic/subsonic regular shock reflection-diffraction configurations, it suffices
to solve the following problem:

\begin{problem}[Free Boundary Problem]\label{fbp-c}
{\it  For $\theta_{\rm w}\in (\theta_{\rm w}^{\rm d}, \frac{\pi}{2})$,
find a free boundary $($curved reflected shock$)$ $\Shock \subset \Lambda\cap \{\cxi<\xi_{1\PtUpL}\}$
{\rm (}$\Shock=P_1P_2$ on {\rm Fig. \ref{fig:RegularReflection}}
and $\Shock=P_0P_2$ on  {\rm Fig.  \ref{fig:SubsonicRegularReflection}}$)$
and a function $\varphi$ defined in the domain
$\Omega$ as shown in  {\rm Figs. \ref{fig:RegularReflection}}--{\rm \ref{fig:SubsonicRegularReflection}}
such that
\begin{itemize}
\item[\rm (i)]
Equation \eqref{2-1} is satisfied in $\Omega$, and the equation is strictly elliptic for $\varphi$ in $\overline\Omega\setminus\overline\Sonic${\rm ,}
\item[\rm (ii)]
$\vphi=\vphi_1$ and $\rho \D\vphi\cdot\nnu_{\rm s}=\rho_1 \D\vphi_1\cdot\nnu_{\rm s}$ {on} the free boundary $\Shock$,
\item[\rm (iii)]
$\vphi=\vphi_2$ and $\D\vphi=\D\vphi_2$ {on} $\Gamma_{\rm sonic}$
in the supersonic case as shown in  {\rm Fig. \ref{fig:RegularReflection}}
 and at $P_0$ in the subsonic case as shown in {\rm Fig. \ref{fig:SubsonicRegularReflection}},
\item[\rm (iv)]
$\D\vphi\cdot\nnu_{\rm w}=0$ {on} $\Wedge=P_0P_3$, and
$\D\vphi\cdot\nnu_{\rm sym}=0$ {on} $\Symm$,
\end{itemize}
where $\nnu_{\rm s}$
is the interior unit normal to $\Omega$
on $\shock$.
}
\end{problem}

Indeed, if $\varphi$ is a solution of Problem \ref{fbp-c}, we extend $\varphi$
from $\Omega$ to $\Lambda$ to become a global entropy solution
(see Figs. {\rm \ref{fig:RegularReflection}}--{\rm \ref{fig:SubsonicRegularReflection}})
by defining:
\begin{equation}\label{phi-states-0-1-2-MainThm}
\varphi=\begin{cases}
\, \varphi_0 \qquad\, \mbox{for}\,\, \xi_1>\xi_1^0 \mbox{ and } \xi_2>\xi_1\tan\theta_{\rm w},\\
\, \varphi_1 \qquad\, \mbox{for}\,\, \xi_1<\xi_1^0
  \mbox{ and above curve} \,\, P_0\PtUpL\PtLwL,\\
\, \varphi_2 \qquad\, \mbox{in region}\,\, P_0\PtUpL\PtUpR.
\end{cases}
\end{equation}
For the subsonic reflection case,
domain $P_0\PtUpL\PtUpR$ is one point
and curve $P_0\PtUpL\PtLwL$ is $P_0\PtLwL$.
Then the global solutions involve two types of transonic (hyperbolic-elliptic) transition:
one is from the {\it hyperbolic} to the {\it elliptic} phases via $\Gamma_{\rm shock}$;
the other is from the {\it hyperbolic}
to the {\it elliptic} phases via $\Gamma_{\rm sonic}$.

The conditions in  {Problem \ref{fbp-c}}(ii)
are the Rankine-Hugoniot conditions \eqref{1-i}--\eqref{1-h} on $\Shock$ between
$\varphi_{|\Omega}$ and $\varphi_1$.
Since $\Shock$ is a free boundary and equation  \eqref{2-1} is strictly elliptic for $\varphi$
in $\overline\Omega\setminus\overline\Sonic$,
then two conditions --- the Dirichlet and oblique derivative conditions --- on $\Shock$ are consistent with
one-phase free boundary problems for nonlinear elliptic PDEs of second order.

In the supersonic case,
the conditions in {Problem \ref{fbp-c}}(iii) are the Rankine-Hugoniot conditions on $\Sonic$ (weak discontinuity)
between
$\varphi_{|\Omega}$ and $\varphi_2$
so that, if $\varphi$ is a solution of {Problem \ref{fbp-c}},
its extension by \eqref{phi-states-0-1-2-MainThm} is a weak solution of {Problem \ref{bvp-c}}.
Since $\Sonic$ is not a free boundary (its location is fixed), it is impossible in general to prescribe
the two conditions given in  {Problem \ref{fbp-c}}(iii) on $\Sonic$ for a second-order elliptic PDE.
In the iteration problem,
we prescribe the condition: $\varphi=\varphi_2$ on $\Sonic$, and then prove that
$\D\varphi=\D\varphi_2$ on $\Sonic$ by exploiting the elliptic degeneracy on $\Sonic$.

\begin{figure}
\centering
\begin{minipage}{0.475\textwidth}
\centering
\includegraphics[height=1.59in,width=1.90in]{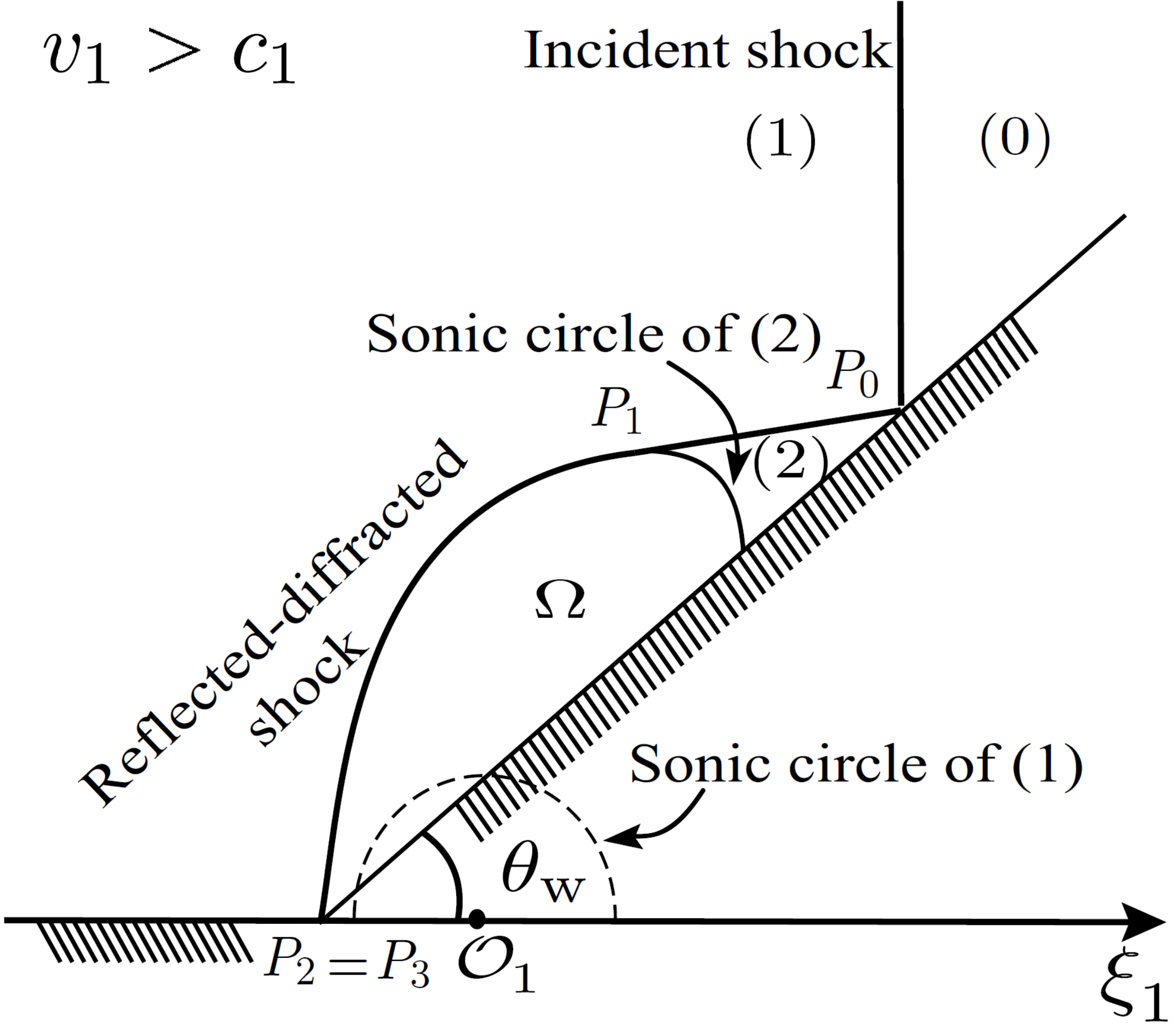}
\caption{\scriptsize Attached supersonic regular reflection-diffraction configuration}
\label{fig:RegularReflection-a}
\end{minipage}
\quad\,\,
\begin{minipage}{0.475\textwidth}
\centering
\includegraphics[height=1.50in,width=1.90in]{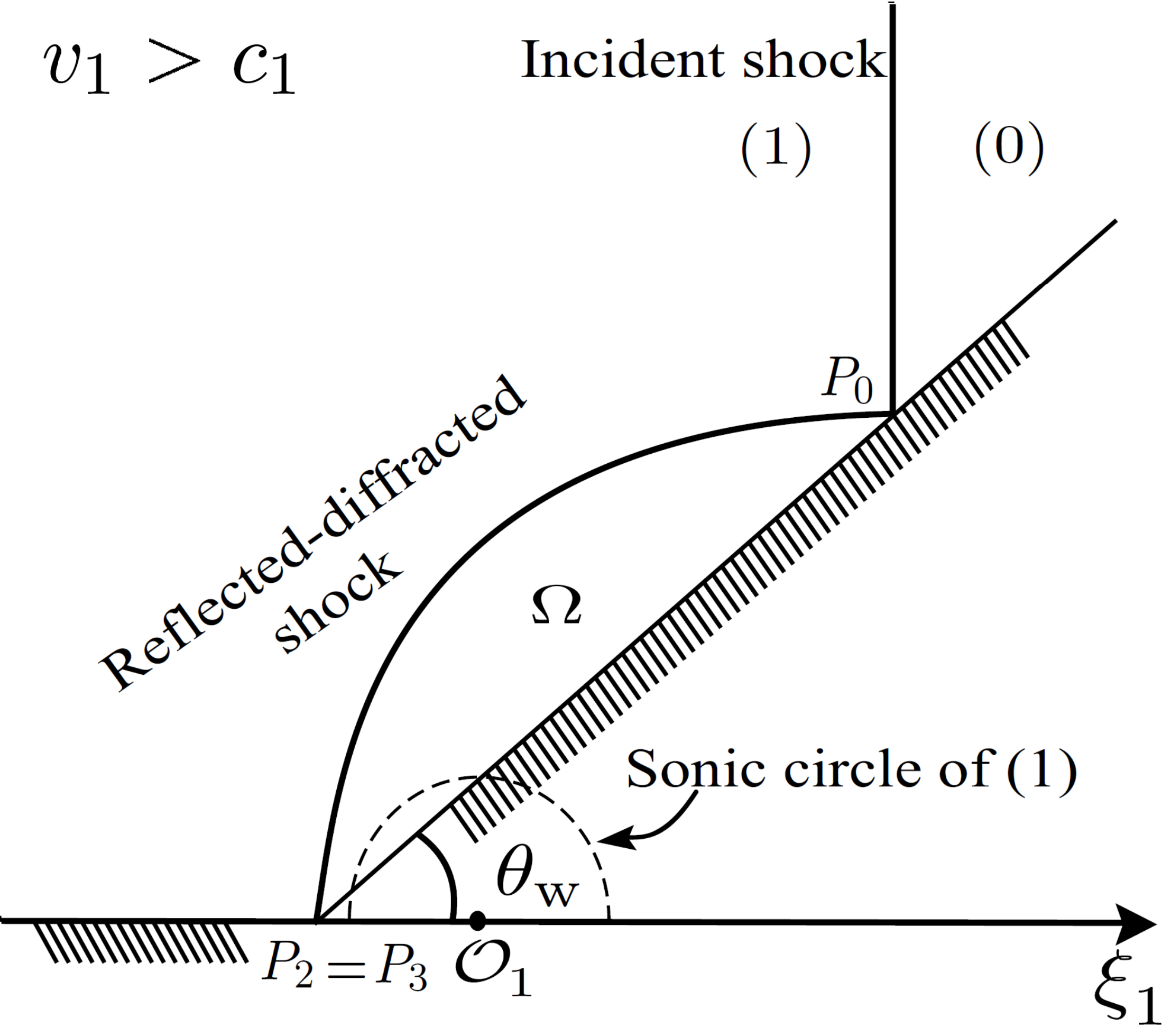}
\caption{\scriptsize Attached subsonic regular reflection-diffraction configuration}\label{fig:SubsonicRegularReflection-a}
\end{minipage}
\vspace{-6pt}
\end{figure}

\smallskip
We observe that there is an additional possibility to the regular shock reflection-diffraction configurations
(beyond the conjectures by von Neumann \cite{vN1}):
for some wedge angle $\theta_{\rm w}^{\rm a}\in (\theta_{\rm w}^{\rm d}, \frac{\pi}2)$,
$\Gamma_{\rm shock}$
may attach to the wedge vertex $\PtLwR$, as observed
by experimental results ({\it cf.} \cite{CF18});
also see Figs. \ref{fig:RegularReflection-a}--\ref{fig:SubsonicRegularReflection-a}.
To describe the conditions of such an attachment,
we use the explicit expressions \eqref{ID-1}
to see that, for each $\rho_0$,
there exists $\rho^{\rm c}>\rho_0$ such that
\begin{eqnarray*}
v_1\le c_1 \quad \mbox{if $\rho_1\in (\rho_0, \rho^{\rm c}]$}; \, \qquad\quad
v_1>c_1 \quad \mbox{if $\rho_1\in (\rho^{\rm c}, \infty)$}.
\end{eqnarray*}
If $v_1\le c_1$, we can rule out the solution with a shock attached to $P_3=(0,0)$.
This is based on the fact that, if $v_1\le c_1$, then
$P_3$ lies within the sonic circle $\overline{B_{c_1}(\mathbf{v}_1)}$ of state (1),
and $\Shock$ does not intersect $\overline{B_{c_1}(\mathbf{v}_1)}$, as we show below.
If $v_1> c_1$, there would be a possibility that
$\Gamma_{\rm shock}$
could be attached to $P_3$
as the experiments show.
With these, the
following results have been obtained:

\smallskip
\begin{theorem}[Chen-Feldman \cite{CF2010,CF18}] \label{mainShockReflThm}
There are two cases{\rm :}
\begin{enumerate}
\item[\rm (i)]
If  $\rho_0$ and $\rho_1$ are such that $v_1\le c_1$, then the
supersonic/subsonic regular reflection solution
exists for each $($half-wedge$)$ angle $\theta_{\rm w}\in (\theta_{\rm w}^{\rm d}, \frac{\pi}{2})$.
That is, for each $\theta_{\rm w}\in (\theta_{\rm w}^{\rm d}, \frac{\pi}{2})$,
there exists a solution $\varphi$ of {\rm Problem \ref{fbp-c}} such that
$$
\Phi(t, {\bf x}) =t\,\varphi(\frac{\bf x}{t}) +\frac{|{\bf x}|^2}{2t}
\qquad\mbox{for}\,\, \frac{\bf x}{t}\in \Lambda,\, t>0,
$$
with
$\rho(t, {\bf x})=\big(\rho_0^{\gamma-1}-(\gamma-1)\big(\partial_t\Phi
      +\frac{1}{2}|\nabla_{\bf x}\Phi|^2\big)\big)^{\frac{1}{\gamma-1}}$,
is a global weak solution of {\rm Problem \ref{ibvp-c}}
satisfying the entropy condition{\rm ;}
that is, $\Phi(t, {\bf x})$ is an entropy solution.

\smallskip
\item[\rm (ii)]
If  $\rho_0$ and $\rho_1$ are such that $v_1> c_1$, then there exists
$\theta_{\rm w}^{\rm a}\in [\theta_{\rm w}^{\rm d}, \frac{\pi}2)$ so that
the regular reflection solution exists for
each angle $\theta_{\rm w}\in (\theta_{\rm w}^{\rm a}, \frac{\pi}{2})$,
and the solution is of the self-similar structure
described in {\rm (i)}  above.
Moreover, if $\theta_{\rm w}^{\rm a}>\theta_{\rm w}^{\rm d}$,
then, for the wedge angle $\theta_{\rm w}=\theta_{\rm w}^{\rm a}$,
there exists an {\it attached} solution, {\it i.e.},  $\varphi$ is
a solution of {\rm Problem \ref{fbp-c}} with  $\PtLwL=\PtLwR$.
\end{enumerate}
The type of regular shock reflection-diffraction
configurations $($supersonic as in {\rm Fig. \ref{fig:RegularReflection}} and {\rm Fig. \ref{fig:RegularReflection-a}},
or  subsonic as in {\rm Fig. \ref{fig:SubsonicRegularReflection}} and {\rm Fig. \ref{fig:SubsonicRegularReflection-a}}$)$
is determined by the type of state
{\rm (2)} at $\PtIncW${\rm :}
\begin{enumerate}[\rm (a)]
\item
For the supersonic/sonic reflection case,
the reflected-diffracted shock $\PtIncW\PtLwL$ is $C^{2,\alpha}$--smooth for some $\alpha\in(0,1)$ and
its curved part $\Gamma_{\rm sonic}$ is $C^\infty$ away from $P_1$.
The solution $\varphi$ is in $C^{1,\alpha}(\overline\Omega)\cap C^\infty(\Omega)$, and is
$C^{1,1}$ across $\Gamma_{\rm sonic}$ which is optimal{\rm ;} that is, $\varphi$ is {\em not} $C^2$ across
$\Gamma_{\rm sonic}$.
\item
For the subsonic reflection case
$($as in {\rm Fig. \ref{fig:SubsonicRegularReflection}} and {\rm Fig. \ref{fig:SubsonicRegularReflection-a}}$)$,
the reflected-diffracted shock $\PtIncW\PtLwL$
and solution $\varphi$ in $\Omega$ are in $C^{1,\alpha}$
near $P_0$ and $P_3$ for some $\alpha\in(0,1)$, and $C^\infty$ away from $\{P_0,P_3\}$.
\end{enumerate}
Moreover, the regular reflection solution tends to the unique normal
reflection $($as in {\rm Fig. \ref{NormReflFigure})}
when the wedge angle $\theta_{\rm w}$ tends to $\frac{\pi}{2}$.
In addition, for both supersonic and subsonic reflection cases,
\begin{align*}
\varphi_2<\varphi<\varphi_1 \quad\mbox{in $\Omega$},
\qquad\quad\,
\D(\varphi_1-\varphi)\cdot {\bf e}\le 0 \quad\mbox{in $\overline\Omega\,\,$ for
  all ${\bf e}\in \overline{Cone({\bf e}_{\xi_2}, {\bf e}_{\CS_1})}$},\label{coneOfMonotRegRefl-cone}
\end{align*}
where $Cone({\bf e}_{\xi_2}, {\bf e}_{\CS_1}):=\{a\,{\bf e}_{\xi_2}+b \,{\bf e}_{\CS_1}\;:\; a, b>0\}$
with ${\bf e}_{\xi_2}=(0,1)$ and ${\bf e}_{\CS_1}$ as the tangent unit vector to $\CS_1$.
\end{theorem}

Theorem \ref{mainShockReflThm} was established by solving {Problem \ref{fbp-c}}.
The first results on the existence of global solutions of the free boundary problem ({Problem \ref{fbp-c}})
were obtained for the wedge angles sufficiently close to $\frac \pi 2$ in Chen-Feldman \cite{CF2010}.
Later, in Chen-Feldman \cite{CF18}, these results
were extended up to the detachment angle
as stated in Theorem \ref{mainShockReflThm}. For this extension, the
techniques developed in  \cite{CF2010}, notably the estimates
near $\Gamma_{\rm sonic}$, were the starting point.

To establish Theorem \ref{mainShockReflThm}, a theory for free boundary problems for nonlinear PDEs of mixed elliptic-hyperbolic type
has been developed, including new methods, techniques, and related ideas.  Some features of these methods and techniques include:

(i)	Exploitation of the global structure of solutions to ensure that the nonlinear PDE \eqref{2-1}
is elliptic for the regular reflection solution in $\Omega$ enclosed by the free boundary $\Gamma_{\rm shock}$ and the fixed boundary
for all wedge angles $\theta_{\rm w}$ up to the detachment angle $\theta_{\rm w}^{\rm d}$ for all the physical cases;
see Figs. 16--19.

(ii) Optimal regularity estimates for the solutions of the {\it degenerate elliptic PDE} \eqref{2-1} both near $\Gamma_{\rm sonic}$
and at corner $P_1$ between the free boundary  $\Gamma_{\rm shock}$
and the elliptic degenerate fixed boundary $\Gamma_{\rm sonic}$ for the supersonic reflection case; see Fig. 16 and Fig. 18.

(iii) For fixed incident shock strength and $\gamma>1$, the dependence of the structural transition of the global solution configurations
on the wedge angle $\theta_{\rm w}$ from the supersonic to subsonic reflection cases, {\it i.e.},
from the degenerate elliptic to the uniformly elliptic
equation \eqref{2-1} near a part of the boundary.

(iv) Uniform {\it a priori} estimates required for all stages of the structural transition between the different configurations.

Based on these methods and techniques for establishing Theorem \ref{mainShockReflThm},
further approaches and related techniques have been developed
to prove that the steady weak oblique transonic shocks, discussed in \S2,
are attainable as large-time asymptotic states by constructing the global  Prandtl-Meyer reflection configurations
in self-similar coordinates in Bae-Chen-Feldman \cite{BCF-14} and the references cited therein,
and that
all the self-similar transonic shocks and related free boundaries in these problems are always convex in Chen-Feldman-Xiang \cite{ChenFeldmanXiang}.

Such questions
also arise in other shock reflection/diffraction problems,
which can be formulated as free boundary problems for transonic shocks
for nonlinear PDEs of mixed type.
These problems
have the following important
attributes: They are physically fundamental and have a wealth of
experimental/numerical data indicating diverse patterns
of complicated configurations
({\it cf.} Figs.
14--\ref{fig:SubsonicRegularReflection-a});
their solutions are building
blocks and asymptotic attractors of general solutions of M-D
hyperbolic conservation laws whose mathematical theory is also
in its infancy ({\it cf.} \cite{BD,CF18,Da,GM}).

\smallskip
Similarly, for the full Euler case, a self-similar solution is a solution of form:
$(\V, p,\rho)(t,\x)=(\vv-\xxi, p,\rho)(\xxi), \, \xxi=\x/t$, governed by
\vspace{-2mm}
\begin{equation}\label{Euler4a}
\begin{cases}
\nabla\cdot(\rho \V)+ n\rho=0,\\
\nabla\cdot(\rho \V\otimes \V) +\nabla p+ (n+1)\rho \V=0,\\
\nabla\cdot\big(\rho \V(E+\frac{p}{\rho})\big)+ n \rho (E+\frac{p}{\rho})=0.
\end{cases}
\end{equation}
System \eqref{Euler4a} is a {\it system of conservation laws of mixed-composite elliptic-hyperbolic type} -- It is
\begin{itemize}
\item strictly {\it hyperbolic} when  $|\V|>c:=\sqrt{\gamma p/\rho}$ (pseudo-supersonic),

\smallskip
\item {\it mixed-composite elliptic-hyperbolic} (two of them are elliptic and the others hyperbolic)
when  $|\V|<c:=\sqrt{\gamma p/\rho}$ (pseudo-subsonic).
\end{itemize}
The transition boundary between
the pseudo-supersonic/pseudo-subsonic phase is $|\V|=c$,
a degenerate set of the solution of equation \eqref{Euler4a}, which is {\it a priori} unknown.

Similar fundamental
mixed problems arise in other applications,
where such nonlinear PDEs of mixed type are the core parts
in even more sophisticated systems; for example, the relativistic Euler equations,
the Euler-Poisson equations, and the Euler-Maxwell equations.

\section{Nonlinear PDEs of Mixed Type and Isometric Embedding Problems in Differential Geometry and Related Areas}

\medskip
Nonlinear PDEs of mixed type also arise naturally from many longstanding
problems in differential geometry
and related areas.
In this section, we first show how the fundamental problem -- the {\it isometric embedding problem} --
in differential geometry can be formulated as
problems for nonlinear PDEs of mixed type, or even no type.

The isometric embedding problem can be stated as follows:
{\it Seek an embedding/immersion of an $n$-D
$($semi-$)$Riemannian
manifold $(\M^n,g)$ with metric $g=(g_{ij})>0$
into an $N$-D
$($semi-$)$
Euclidean space so that the metric, often along with assigned regularity/curvatures, is preserved}.

\smallskip
This problem has assumed a position of fundamental conceptual importance in differential geometry
since the works of Darboux (1894), Weyl (1916), Janet (1926), and Cartan (1927).
A classical question is whether a smooth Riemannian manifold $(\M^n,g)$ can be embedded into $\R^N$
with sufficiently large $N$; {\it e.g.}  Nash (1956), Gromov (1986), and G\"{u}nther (1989).
A further fundamental issue is whether $(\M^n,g)$ can be embedded/immersed in $\R^{s_n}$ with critical Janet
dimension $s_n=\frac{n(n + 1)}{2}$ and assigned regularity/curvatures.
The solution to this issue will
advance our understanding of properties
of (semi-)Riemannian manifolds and provide frameworks/approaches for real applications, including
the problems for realization/stability/rigidity/classification of isometric embeddings
in many important application areas
(see {\it e.g.}
elasticity, materials science, optimal designs,
thin shell/biological leaf growth, protein folding, cell/tissue organization, and manifold data analysis).

\begin{figure}
\begin{minipage}{0.480\textwidth}
\vspace{-10pt}
\centering
\includegraphics[width=3.2cm,height=1.70cm]{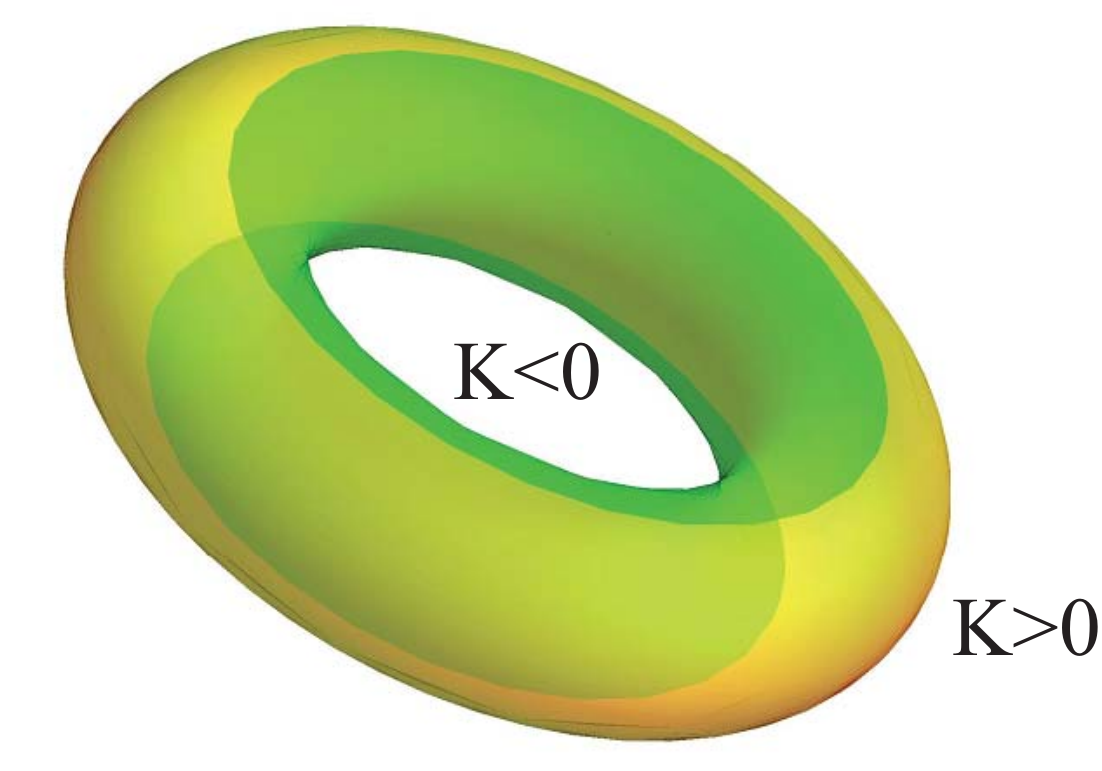}
\caption{\scriptsize \,
The Gauss curvature $K$ of a torus with mixed sign}
\label{torus}
\end{minipage}
\begin{minipage}{0.480\textwidth}
\centering
\includegraphics[height=1.60in,width=1.25in]{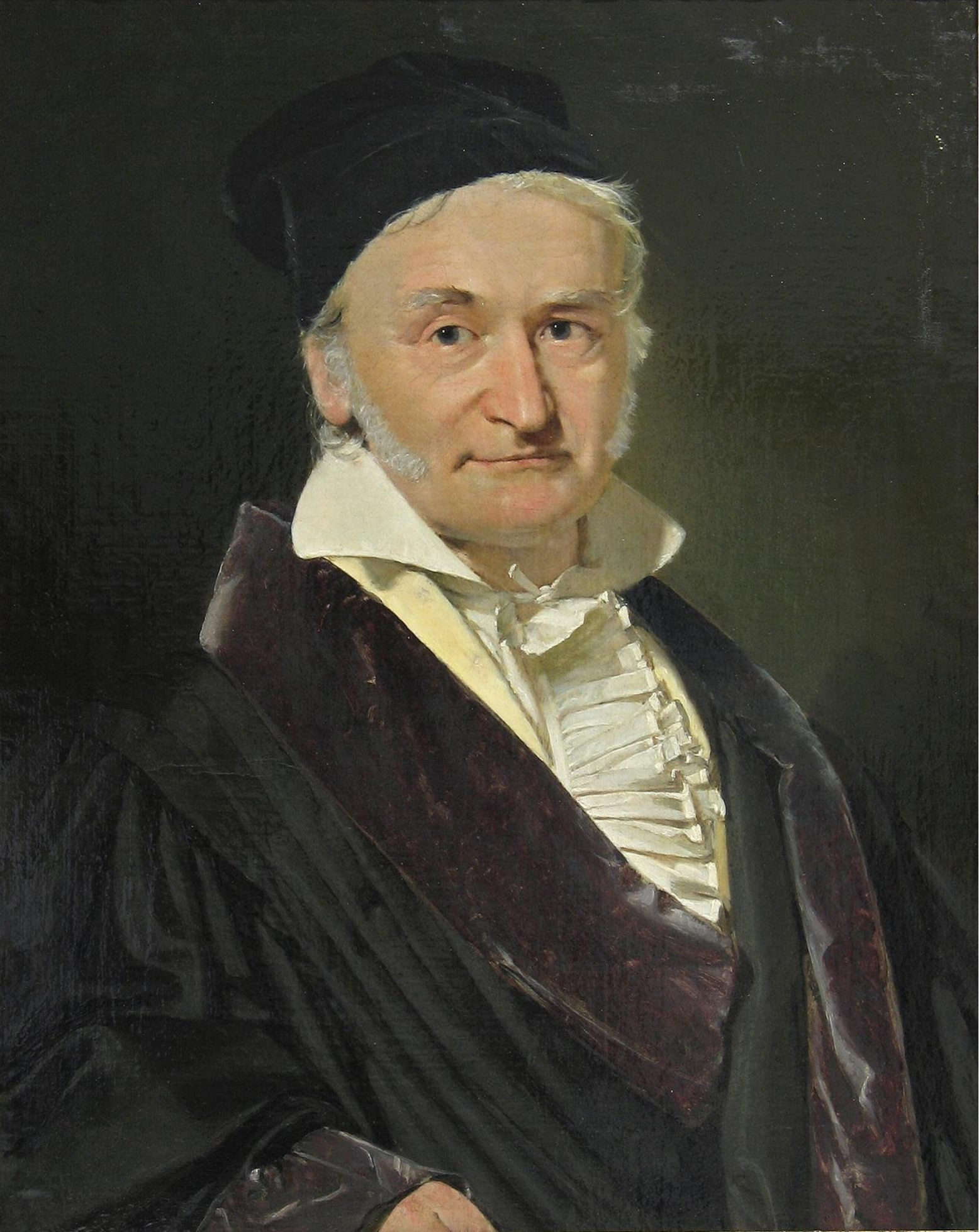}
\caption{\scriptsize Johann Carl Friedrich Gauss (30 April 1777 -- 23 February 1855)
introduced the notion of Gauss (or Gaussian) curvature and the {\it Theorema Egregium}}
\label{Euler}
\end{minipage}
\end{figure}

When $n=2$, following Darboux\footnote{Darboux, G.:
\newblock {\em Le{\c c}ons sur la Th\'{e}orie des Surface},
Vol. {\bf 3}, Gauthier-Villars: Paris, 1894.},
the isometric embedding problem on a chart can be reduced to finding a function $u$
that solves the following nonlinear Monge-Amp\`{e}re equation ({\it cf}. \cite{HH}):
\begin{equation}\label{Monge-eq}
\qquad\quad  \det (\nabla^2 u) = |g|\big(1 -|\nabla u|^2_g\big)K,
\end{equation}
with $|g|=\det(g)$,
$|\nabla u|_g:=\frac{1}{|g|}
  \big(g_{22}|\partial_{x_1}u|^2-2g_{12}\partial_{x_1}\partial_{x_2}u+g_{11}|\partial_{x_2}u|^2\big)<1$
as required, and the Gauss curvature $K=K(g)$ of metric $g$.
Equation \eqref{Monge-eq} is {\it elliptic} if $K>0$,
{\it hyperbolic} if
$K < 0$, and {\it degenerate} when $K=0$.
The sign change of $K$ is very common for surfaces
and is necessary for many important cases;
the simplest example of such surfaces is the torus as shown in Fig. \ref{torus}.

\begin{figure}
\begin{minipage}{0.495\textwidth}
\vspace{-4pt}
\centering
\includegraphics[height=1.50in,width=1.1in]{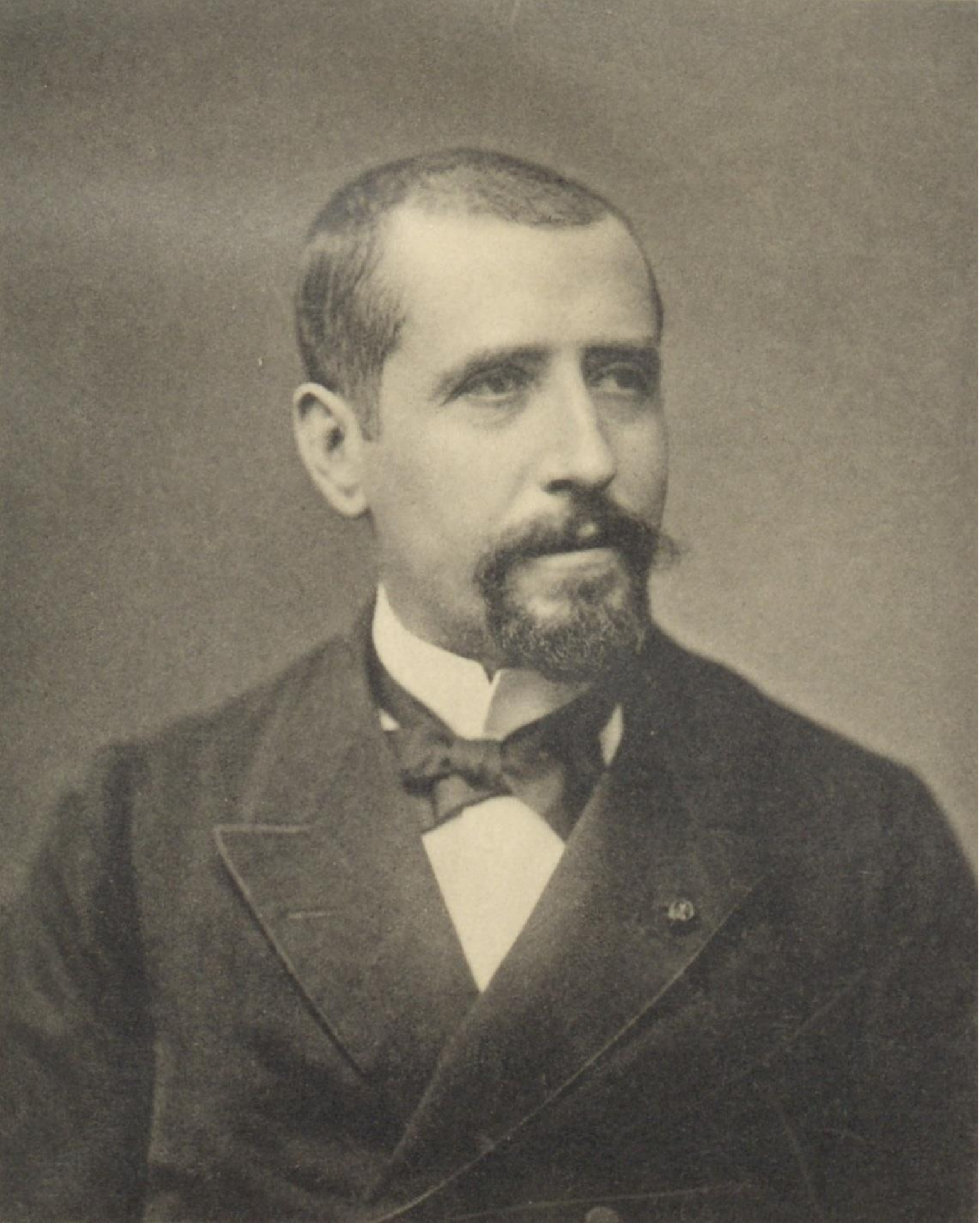}
\caption{\scriptsize \,\,
Jean-Gaston Darboux (14 August 1842 -- 23 February 1917)
indicated the connection between the isometric embedding and
the nonlinear Monge-Amp\`{e}re equation}
\label{Darboux}
\end{minipage}
\begin{minipage}{0.480\textwidth}
\vspace{-18pt}
\centering
\includegraphics[height=1.50in,width=1.1in]{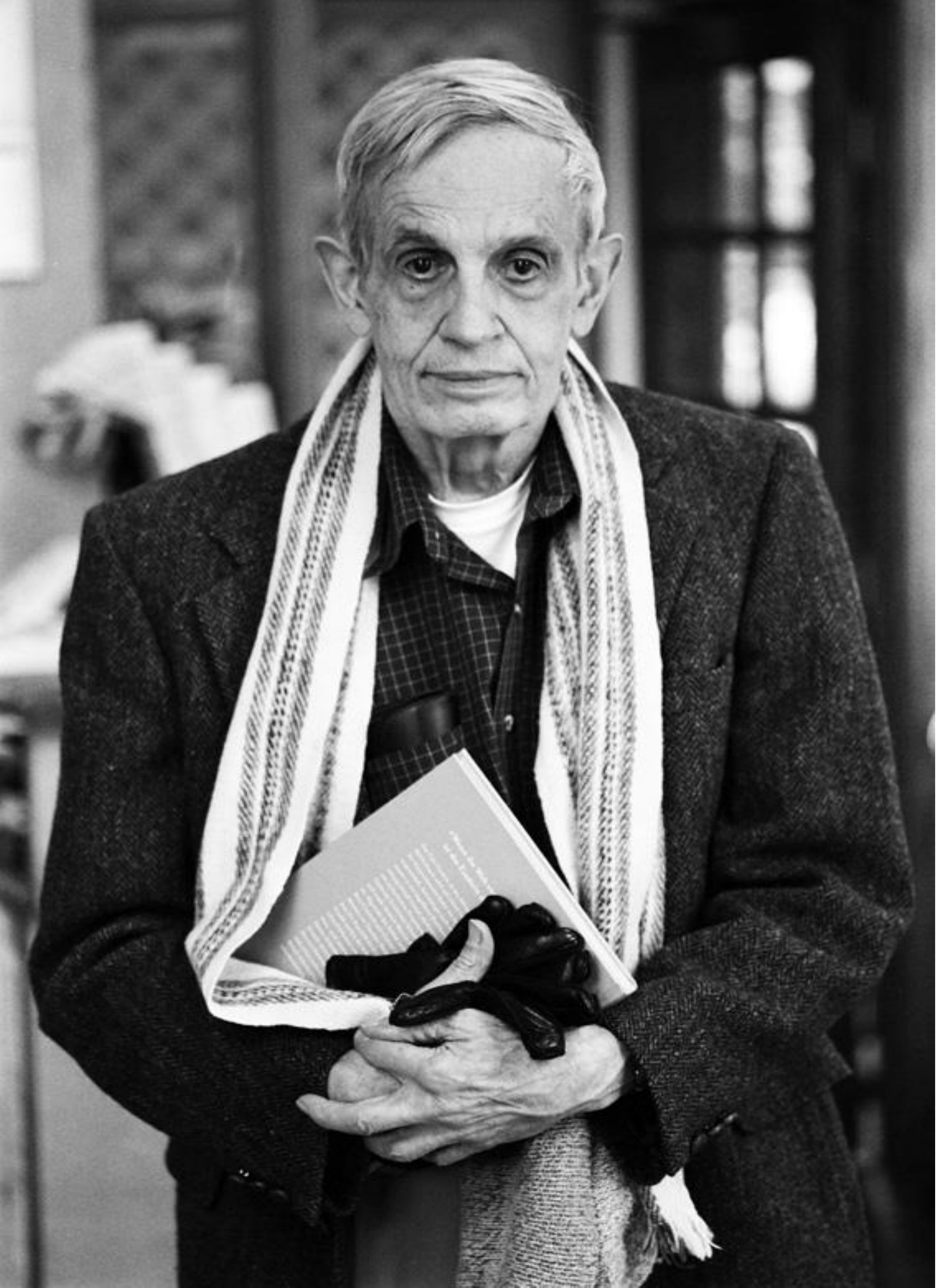}
\caption{\scriptsize John Forbes Nash Jr. (13 June 1928 -- 23 May 2015)
established the Nash embedding theorems}
\label{Nash}
\end{minipage}
\end{figure}

Nirenberg (1953) first solved the {\it Weyl problem},
establishing that
any smooth metric $g$ on $\mathbb{S}^2$ can be
globally
embedded into $\R^3$ smoothly if the Gauss curvature $K>0$.
Then
one would ask whether any 2-D Riemannian surface is always
embeddable into ${\mathbb{R}}^3$.
The answer is {\it no} if $K\leq 0$.
The embedding problem is still largely open for global results for general $K$,
even though some local results have been obtained;
see \cite{HH} and the references therein.

On the other hand, the fundamental theorem of surface theory states that
{\it there exists a simply connected surface in $\R^3$ whose first and second
fundamental forms are $I=g_{ij}dx_idx_j$ and $I\!I=h_{ij}dx_i dx_j$ on a
domain
for $i,j=1,2$,
provided that the coefficients
 $\{h_{ij}\}$, together with metric $g=(g_{ij})>0$,
 satisfy the Gauss-Codazzi equations}:
\vspace{-3mm}
\begin{eqnarray}
&&\,\,LN-M^2=K, \label{g2}\\[0.5mm]
&&\begin{cases}\label{g1}
\d_{x_1}{N}-\d_{x_2}{M}=-\G^{1}_{22}L +2\G^{1}_{12}M -\G^{1}_{11}N,\\
\d_{x_1}{M}-\d_{x_2}{L}=\G^{2}_{22}L-2\G^{2}_{12}M+\G^{2}_{11}N,
\end{cases}
\end{eqnarray}
where $L={h_{11}}/{\sqrt{|g|}}, M={h_{12}}/{\sqrt{|g|}}$, and $N={h_{22}}/{\sqrt{|g|}}$,
and $\G^{k}_{ij}$  are the Christoffel symbols for $i,j,k=1,2$.
This theorem still holds for immersion\footnote{Mardare, S.: The fundamental
theorem of surface theory for surfaces with little regularity,
{\em J. Elasticity}, {\bf 73} (2003), 251--290.}
even when $\{h_{ij}\}$ is only in
$L^p$ for $p>2$.
Thus, given
$(g_{ij})>0$, system \eqref{g2}--\eqref{g1}
consists of three nonlinear PDEs for the unknowns $(L, M, N)$
determining $\{h_{ij}\}$ whose knowledge
gives a desired immersion.
Then the problem can be reduced to the solvability
of system
\eqref{g1} under constraint \eqref{g2}, which is of {\it mixed elliptic-hyperbolic type} determined by
the sign of the Gauss curvature $K$.
From the viewpoint of
geometry, \eqref{g2} is a constraint condition,
while \eqref{g1}
are compatibility conditions.

System \eqref{g2}--\eqref{g1} has similar features to those
in gas dynamics in \S 2--\S 3.
A natural question is whether it can be written in a gas dynamic formulation
to examine underlying interrelated connections.
Indeed,
a novel observation in Chen-Slemrod-Wang \cite{CSW10}
has indicated that
this is
the case: the Codazzi system \eqref{g1} can be formulated as
the familiar nonlinear balance laws of momentum:
\begin{equation} \label{g3}
\vspace{-3pt}
\begin{cases}
\d_{x_1}(\r u^2+p)+\d_{x_2}(\r uv)
 =-\G^{1}_{22}(\r v^2+p)-2\G^{1}_{12}\r uv-\G^{1}_{11}(\r u^2+p),\\[1.5mm]
 \d_{x_1}(\r uv)+\d_{x_2}(\r v^2+p)
 =-\G^{2}_{22}(\r v^2+p)-2\G^{2}_{12}\r uv-\G^{2}_{11}(\r u^2+p),
\end{cases}
\end{equation}
and the Gauss equation \eqref{g2} becomes the {\it Bernoulli} relation:
$
\r=(q^2+K)^{-1/2}
$
if $p=-1/\r$ chosen as
the Chaplygin pressure for $q=\sqrt{u^2+v^2}$.
In this case,  define the sound speed:
$c=\sqrt{p'(\r)}=1/\r$.
Then
\begin{itemize}
\item $q<c$ and the {\it flow} is subsonic  when $K>0$,
\item $q>c$ and the {\it flow} is supersonic  when $K<0$,
\item $q=c$ and the {\it flow} is sonic when $K=0$.
\end{itemize}

A weak compactness framework has been introduced and applied for establishing
the existence and weak continuity/stability of isometric embeddings in $W^{2,p}, p\ge 2$, in \cite{CL17,CSW10},
which has shown its high potential.
In particular, the weak continuity/stability of the Gauss-Codazzi equations \eqref{g2}--\eqref{g1}
and isometric immersions of (semi-)Riemannian manifolds, independent of local coordinates,
have been established in \cite{ChenGiron,CL17}, even for the case
$p=2$.

\smallskip
For the higher dimensional case, the Gauss-Codazzi equations for ${\boldsymbol h}=\{h_{ij}^{a}\}$ are coupled with
the Ricci equations
for
the coefficients ${\boldsymbol \kappa}=\{\kappa_{lb}^{a}\}$ of
the connection form on
the normal bundle to become
the Gauss-Codazzi-Ricci equations in a local coordinate chart of the manifold:
\begin{align}
& h_{ji}^a h_{kl}^a-h_{ki}^a h_{jl}^a=R_{ijkl} \quad\text{(Gauss equations)}\label{G1},\\
&\d_{x_k} h_{lj}^a- \d_{x_l} h_{kj}^a
=-\Gamma_{lj}^mh_{km}^a+\Gamma_{kj}^m h_{lm}^a -\big(\kappa_{kb}^a h_{lj}^b-\kappa_{lb}^ah_{kj}^b\big)
  \quad \text{(Codazzi equations)}, \label{C1}\\
&\d_{x_k}\kappa_{lb}^a- \d_{x_l}\kappa_{kb}^a
= g^{ij}\big(h^a_{li}h^b_{kj}-h^a_{ki}h^b_{lj}\big)
  +\kappa_{lc}^a\kappa_{kb}^c-\kappa_{kc}^a\kappa_{lb}^c\quad \text{(Ricci equations)},\label{R1}
\end{align}
where
$\kappa_{kb}^a=-\kappa_{ka}^b$ are the coefficients of the
connection form
on the normal bundle,
$R_{ijkl}$ is the Riemann curvature tensor,
the
indices $a, b, c$ run from $1$ to $N$, and $i, j, k, l, m, n$
run from $1$ to $d\ge 3$.
System \eqref{G1}--\eqref{R1} has no
type, neither purely hyperbolic nor purely elliptic
for general Riemann curvature tensor $R_{ijkl}$.
Nevertheless, the weak continuity of the nonlinear
system \eqref{G1}--\eqref{R1}
has been established.

\begin{theorem}[Chen-Slemrod-Wang \cite{CSW10}]\label{thm-4.1}
Let $({\boldsymbol h}^\e, {\boldsymbol \kappa}^\e)$
be a
sequence of solutions of the Gauss-Codazzi-Ricci system \eqref{G1}--\eqref{R1}, which is
uniformly bounded in $L^p$ for $p>2$.
Then the weak limit vector
field $({\boldsymbol h}, {\boldsymbol \kappa})$
of the sequence $({\boldsymbol h}^\e, {\boldsymbol \kappa}^\e)$
in $L^p$ is
still a solution to the Gauss-Codazzi-Ricci system \eqref{G1}--\eqref{R1}.
\end{theorem}

The proof is based on our key observation in \cite{CSW10} for
the div-curl structure of system \eqref{G1}--\eqref{R1}:
For fixed $i, j, k, l, a,b, c$,
\vspace{-2mm}
\begin{align}
&\text{div}\,(\underbrace{\overbrace{0, \cdots, 0,
h_{li}^{a,\e}}^{k}, 0, \cdots, -h_{ki}^{a,\e}}_{l},0,\cdots,0)=R_1^{\e},
\quad \text{curl}\,(h^{b,\e}_{1j}, h^{b,\e}_{2j}, \cdots,
h^{b,\e}_{dj})=R_2^\e,\label{R-1}\\[-2mm]
&\text{div}\,(\underbrace{\overbrace{0, \cdots, 0,
\kappa_{lc}^{a,\e}}^{k}, 0, \cdots, -\kappa_{kc}^{a,\e}}_l, 0,
\cdots, 0)=R_3^{\e},
\quad \text{curl}\,(\k^{c,\e}_{1b},\k^{c,\e}_{2b},\cdots,\k^{c,\e}_{db})=R_4^\e, \label{R-2}\\[-2mm]
&\text{div}\,(\underbrace{\overbrace{0, \cdots, 0,
h_{lj}^{b,\e}}^{k}, 0, \cdots, -h_{kj}^{b,\e}}_{l},0,\cdots, 0)=R_5^\e,
\quad
\text{curl}\,(\k^{a,\e}_{1b},
\k^{a,\e}_{2b}, \cdots, \k^{a,\e}_{db})=R_6^{\e}, \label{R-3}
\end{align}
where $R_r, r=1,\dots,6$, consist of the three types of nonlinear
quadratic terms:
$$
{h_{li}^{a,\e}h_{kj}^{b,\e}-h_{ki}^{a,\e}h_{lj}^{b,\e}},
\quad
\kappa_{lc}^{a,\e}\kappa_{kb}^{c,\e}-\kappa_{kc}^{a,\e}\kappa_{lb}^{c,\e}  ,
\quad
{\kappa_{kb}^{a,\e}h_{lj}^{b,\e}-\kappa_{lb}^{a,\e}h_{kj}^{b,\e}},
$$
besides several linear terms involving  $({\boldsymbol h}^\e, {\boldsymbol \kappa}^\e)$,
while the nonlinear quadratic terms are actually the scalar products of the vector fields
given on the left-hand sides of \eqref{R-1}--\eqref{R-3}.
Therefore,
this div-curl structure fits the following classical
div-curl lemma divinely (Murat 1978, Tartar 1979):
{\it Let $\Omega\subset\R^d, d\ge 2,$ be open and bounded. Let $p, q>1$ such
that $\frac{1}{p}+\frac{1}{q}=1$.
Assume that, for $\varepsilon>0$, two fields
$\mathbf{u}^\varepsilon\in L^p(\Omega; \R^d)$
and $\mathbf{v}^\varepsilon\in L^q(\Omega; \R^d)$
satisfy the following conditions{\rm :}

\begin{enumerate}\renewcommand{\theenumi}{\roman{enumi}}
\item[\rm (i)]
\label{weakSol-def-i1}
$\mathbf{u}^\varepsilon\rightharpoonup \mathbf{u}$ weakly in
$L^p(\Omega;\R^d)$
and $\mathbf{v}^\varepsilon\rightharpoonup \mathbf{v}$ weakly in
$L^q(\Omega;\R^d)$ as $\varepsilon\to 0$,

\item[\rm (ii)]  ${\rm div}\, \mathbf{u}^\varepsilon$ are confined in a
compact subset of $W^{-1, p}_{\rm loc}(\Omega; \R)$,

\item[\rm (iii)] ${\rm curl}\, \mathbf{v}^\varepsilon$ are confined in a
compact subset of $W^{-1, q}_{\rm loc}(\Omega; \R^{d\times d})$,
\end{enumerate}
where $W^{-1,p}(\Omega;\R)$ is the dual space of $W^{1,q}(\Omega;\R)$, and vice versa.
Then the scalar product of $\mathbf{u}^\e$ and $\mathbf{v}^\e$ are weakly continuous{\rm :}
$\mathbf{u}^\varepsilon\cdot \mathbf{v}^\varepsilon\longrightarrow \mathbf{u}\cdot \mathbf{v}$
in the sense of distributions.
}

\smallskip
With this div-curl lemma, the weak continuity result in Theorem \ref{thm-4.1}
can be seen as follows: For the uniformly bounded
sequence  $({\boldsymbol h}^\e, {\boldsymbol \kappa}^\e)$
in $L^p, p>2$,
$R_r^\varepsilon, r=1,\dots,6,$ are uniformly bounded in $L^{p/2}$, which implies
that $R_r^\varepsilon, r=1,\dots,6,$ are compact in $W^{-1, q}_{\rm loc}$ for some $q\in (1,2)$.
On the other hand, equations \eqref{R-1}--\eqref{R-3} imply that
$R_r^\varepsilon, r=1,\dots,6,$ are uniformly bounded in $W^{-1, p}_{\rm loc}$ for $p>2$.
Then the interpolation compactness argument can yield that
$$
R_r^\varepsilon, \,\, r=1,\dots,6,
\qquad
\text{are confined in a compact set in $H^{-1}_{\rm loc}(\Omega)$}.
$$
With this, we can employ the div-curl lemma to conclude that
\begin{eqnarray*}
&&({h_{li}^{a,\e}h_{kj}^{b,\e}-h_{ki}^{a,\e}h_{lj}^{b,\e}},\,
\kappa_{lc}^{a,\e}\kappa_{kb}^{c,\e}-\kappa_{kc}^{a,\e}\kappa_{lb}^{c,\e},\,
{\kappa_{kb}^{a,\e}h_{lj}^{b,\e}-\kappa_{lb}^{a,\e} h_{kj}^{b,\e}})\\
&&\longrightharpoonup
({h_{li}^{a}h_{kj}^{b}-h_{ki}^{a}h_{lj}^{b}},\,
{\kappa_{lc}^a\kappa_{kb}^c-\kappa_{kc}^a\kappa_{lb}^c},\,
{\kappa_{kb}^ah_{lj}^b-\kappa_{lb}^ah_{kj}^b})
\end{eqnarray*}
in the sense
of distributions as {$\varepsilon\to 0$}.  Then Theorem \ref{thm-4.1} follows.

This local weak continuity result can be extended to the global weak continuity
of the Gauss-Codazzi-Ricci system \eqref{G1}--\eqref{R1}:

\begin{theorem}[Chen-Li \cite{CL17}]\label{thm:4.4}
Let $(M,g)$ be a Riemannian manifold with $g\in W^{1,p}$ for $p>2$.
Let $({\boldsymbol h}^\e, {\boldsymbol \kappa}^\e)$ be a sequence of solutions
$($i.e., the coefficients of the second
fundamental form and the connection form on the normal bundle$)$
in $L^p$
of the Gauss-Codazzi-Ricci system \eqref{G1}--\eqref{R1} in the distributional sense.
Assume that, for any submanifold $K\Subset M$, there exists $C_K>0$ independent
of $\e$ such that
$$
\sup_{\epsilon>0} \|({\boldsymbol h}^\e, {\boldsymbol \kappa}^{\e})\|_{L^p(K)} \leq C_K.
$$
Then, when $\e\to 0$, there exists a subsequence
of $({\boldsymbol h}^\e, {\boldsymbol \kappa}^\e)$ that
converges weakly in $L^p$ to
a pair $({\boldsymbol h}, {\boldsymbol \kappa})$ that is still a weak
solution of the Gauss-Codazzi-Ricci system \eqref{G1}--\eqref{R1}.
\end{theorem}

The proof is based on a compensated compactness theorem in Banach spaces,
which leads directly to a globally intrinsic div-curl lemma on Riemannian manifolds,
developed in Chen-Li \cite{CL17}.
From the viewpoint of geometry, the $L^p$ bounded requirement on the connection form
on the normal bundle ${\boldsymbol \kappa}^{\e}$ is not intrinsic.
Therefore, Theorem \ref{thm:4.4} has been reformulated in the following theorem:

\begin{theorem}[Chen-Giron \cite{ChenGiron}]\label{thm:4.5}
Let {$(M,g)$} be a Riemannian manifold with {$g\in W^{1,p}$} for {$p>2$}.
Let $({\boldsymbol h}^\e, {\boldsymbol \kappa}^\e)$ be a sequence of solutions
$($i.e., the coefficients of the second
fundamental form and the connection form on the normal bundle$)$ in $L^p$
of the Gauss-Codazzi-Ricci system \eqref{G1}--\eqref{R1} in the distributional sense.
Assume that, for any submanifold $K\Subset M$, there exists $C_K>0$ independent
of $\e$ such that
$$
\sup_{\epsilon>0} \|{\boldsymbol h}^\e\|_{L^p(K)} \leq C_K.
$$
Then there exists a refined sequence
$(\tilde{\boldsymbol h}^\e, \tilde{\boldsymbol \kappa}^\e)$, each of which is still a weak
solution of the Gauss-Codazzi-Ricci system \eqref{G1}--\eqref{R1}, such that,
when $\e\to 0$,
$(\tilde{\boldsymbol h}^\e, \tilde{\boldsymbol \kappa}^\e)$ converges weakly in $L^p$ to
a pair $({\boldsymbol h}, {\boldsymbol \kappa})$ that is still
a weak solution of the Gauss-Codazzi-Ricci system \eqref{G1}--\eqref{R1}.
\end{theorem}

As a direct corollary, the weak limit of isometrically immersed surfaces with lower regularity
in $W^{2,p}$ is still an
isometrically immersed surface in $\R^d$ governed by the Gauss-Codazzi-Ricci system \eqref{G1}--\eqref{R1}
for any $R_{ijkl}$ (without sign/type restriction)
with respect only to the coefficients of the second
fundamental form.
The weak continuity result in Theorem \ref{thm:4.5}
is global and intrinsic, independent of local coordinates, with no restriction on the Riemann curvatures
and the types of system
\eqref{G1}--\eqref{R1}.
The key to the proof is to exploit the invariance for a choice of suitable gauge to control the
full connection form and to develop a non-abelian div-curl lemma on Riemannian manifolds
(see Chen-Giron \cite{ChenGiron}).

This approach and related observations have been
motivated by
the theory of polyconvexity in nonlinear elasticity\footnote{
Ball, J.:
Convexity conditions and existence theorems in nonlinear elasticity,
{\em Arch. Ration. Mech. Anal.} {\bf 63} (1976), 337--403.},
intrinsic methods in elasticity and nonlinear Korn inequalities\footnote{
see Ciarlet, P.~G.:
\textit{Mathematical Elasticity}, Volume 1: \textit{Three--Dimensional Elasticity}.
North-Holland: Amsterdam, 1988;
\textit{An Introduction to Differential Geometry with Applications to Elasticity},
Springer: Dordrecht, 2005.
},
Uhlenbeck compactness
and Gauge theory\footnote{Uhlenbeck, K.~K.: Connections with {$L^{p}$} bounds on curvature,
 {\it Comm. Math. Phys.} {\bf 83} (1982), pp.~31--42.}$^{,}$\footnote{Donaldson, S. K.:
 An application of gauge theory
 to four-dimensional topology,
{\it J. Diff. Geom.} {\bf 18} (1983), 279--315.},
among other ideas.

\section{Further Connections, Unified Approaches, and Current Trends}

In \S 2--\S 4, we have presented several important sets of
nonlinear PDEs of mixed elliptic-hyperbolic type, or even no type,
in shock wave problems in fluid mechanics and isometric embedding problems in
differential geometry and related areas.
Such nonlinear PDEs of mixed type arise naturally in other problems in fluid mechanics,
differential geometry/topology, nonlinear elasticity, materials science, mathematical physics, dynamical systems,
and related areas.

We have shown that free boundary methods, weak convergence methods, and related techniques
are useful as unified approaches to deal with the nonlinear mixed problems
involving {\it both elliptic and hyperbolic phases} in \S 2--\S4.
Friedrichs's positive symmetric techniques have also shown high potential in solving
mixed-type problems\footnote{Chen, G.-Q., Clelland, J., Slemrod, M., Wang, D., Yang, D.:
Isometric embedding via strongly symmetric positive systems, {\it Asian J. Math.} {\bf 22} (2018), 1--40.}.
Entropy methods and kinetic methods have been useful for solving
nonlinear PDEs of hyperbolic or mixed hyperbolic-parabolic type.
Variational approaches deserve to be further explored, especially for handling transonic flow problems since
the solutions of these problems are critical points of the corresponding functionals.
Some approximate methods such as viscosity methods, relaxation methods, shock capturing methods,
stochastic methods, and related numerical methods
should be further analyzed/developed,
and
numerical calculations/simulations should also be performed to gain new ideas and motivations.
These methods, along with
energy estimate techniques,
functional analytical methods, measure-theoretic techniques ({\it esp}. divergence-measure fields),
and other methods, should be further developed into
more powerful approaches, applicable to wider classes of nonlinear PDEs of mixed type.
The underlying structures of the nonlinear PDEs of mixed type
under consideration
have been one of our motivating factors in developing
new methods/techniques/ideas for unified approaches.
As indicated earlier, the analysis of nonlinear PDEs of mixed type is still
in its early stages and most nonlinear mixed-type problems are widely open,
requiring additional new ideas, methods, and techniques.

\bigskip
\noindent
{\bf Acknowledgements}.
The author would like to thank his collaborators and former students, including Myoungjean Bae,
Jun Chen,
Jeanne Clelland,
Mikhail Feldman, Tristin Giron, Siran Li, Marshall Slemrod,
Dehua Wang, Wei Xiang,
and Deane Yang,
as well as the colleagues whose work should be cited (but could not be done so as he wished, due
to the strict limitation of the reference number and the length by the journal; please see the references cited in [1]--[20]),
for their explicit and implicit contributions to the material
presented in this article.
The research of Gui-Qiang G. Chen was supported in part by
the UK Engineering and Physical Sciences Research Council Award
EP/L015811/1, EP/V008854, and EP/V051121/1.


\begin{thebibliography}{10}

\bibitem{BCF-14}
Bae, M., Chen, G.-Q., Feldman, M.:
\textit{Prandtl-Meyer Reflection Configurations,
Transonic Shocks, and Free Boundary Problems},
Research Monograph, 118 pages,
Memoirs of Amer. Math. Soc.,
AMS: Providence, RI, 2022 (to appear).


\bibitem{BD}
Ben-Dor, G.:
\newblock{\itshape Shock Wave Reflection Phenomena},
\newblock{Springer-Verlag: New York}, 1991.

\bibitem{Bers}
Bers, L.:
{\em Mathematical Aspects of Subsonic and Transonic
Gas Dynamics},
John Wiley \& Sons, Inc.: New York, 1958.


\bibitem{Bet}
Bitsadze, A.~V.:
\textit{Equations of the Mixed Type}, Macmillan
Company: New York, 1964.


\bibitem{CF2010}
Chen, G.-Q., Feldman, M.:
Global solutions to shock reflection by large-angle wedges for potential flow,
{\itshape Ann. of Math.} \textbf{171} (2010), 1019--1134.


\bibitem{CF18}
Chen, G.-Q., Feldman, M.:
\newblock {\it Mathematics of Shock Reflection-Diffraction and von Neumann's Conjecture}.
\newblock Research Monograph, Annals of Mathematics Studies, \textbf{197},
Princeton University Press, Princeton, 2018.

\bibitem{CF22}
Chen, G.-Q., Feldman, M.:
\newblock Multidimensional transonic shock waves and free boundary problems,
{\it Bulletin of Mathematical Sciences} \textbf{12}(01), 2230002




\bibitem{ChenFeldmanXiang}
Chen, G.-Q., Feldman, M., Xiang, W.:
Convexity of self-similar transonic shock waves for potential flow,
\textit{Arch. Ration. Mech. Anal.} \textbf{238} (2020), 47--124.


\bibitem{ChenGiron}
Chen, G.-Q.,  Giron, T.:
\newblock
Weak continuity of Gauss-Codazzi-Ricci equations with $L^p$-bounded second fundamental form,
{\em Preprint}, 2022 (to be posted at arXiv).


\bibitem{CL17}
Chen, G.-Q., Li, S.:
\newblock Global weak rigidity of the Gauss-Codazzi-Ricci equations and isometric immersions
of Riemannian manifolds with lower regularity,
{\em J. Geom. Anal.} {\bf 28(3)} (2018), 1957--2007.



\bibitem{CSW10}
Chen, G.-Q., Slemrod, M., Wang, D.-H.:
\newblock
Isometric embedding and compensated compactness,
{\em CMP.} {\bf 294} (2010), 411--437; Weak continuity of the Gauss-Codazzi-Ricci
system for isometric embedding,
{\em Proc. Amer. Math. Soc.} {\bf 138} (2010), 1843--1852.

\bibitem{CFr}
Courant, R., Friedrichs, K.~O.:
{\itshape Supersonic Flow and Shock Waves},
Springer-Verlag: New York, 1948.


\bibitem{Da} Dafermos, C. M.:
\newblock{\it Hyperbolic Conservation Laws in Continuum Physics},
4th Ed., Springer-Verlag: Berlin, 2016.


\bibitem{Evans}
Evans, L.~C.:
{\it  Partial Differential Equations}.
Second edition. Graduate Studies in Mathematics, 19.
American Mathematical Society: Providence, RI, 2010.


\bibitem{GM}
Glimm, J., Majda, A.:
\newblock
{\em Multidimensional Hyperbolic Problems and Computations},
Springer-Verlag: NY, 1991.


\bibitem{Had}
Hadamard, J.:
{\it Lectures on Cauchy's Problem in Linear Partial Differential Equations}.
Yale University Press: New York,
 1923;
Dover Publications, New York, 1953.


\bibitem{HH}
Han, Q. and Hong, J.-X.:
\newblock {\em Isometric Embedding of Riemannian Manifolds in Euclidean
Spaces}, AMS, RI, 2006.


\bibitem{Mo85}
Morawetz, C.~S.:
{On a weak solution for a transonic flow problem},
\textit{Comm. Pure Appl. Math.} \textbf{38}: 797--818, 1985.

\bibitem{vN1}
 von Neumann, J.:
{Oblique reflection of shocks}, {\itshape Explo. Res. Rep.} \textbf{12}, Navy
Department, Bureau of Ordnance, Washington, DC, 1943;
Refraction, intersection, and reflection of shock waves,
{\em NAVORD Rep.} \textbf{203-45}, Navy Department, Bureau of Ordnance,
Washington, DC, 1945;
{\em Collected Works}, Vol 6, Pergamon Press, 1963.

\bibitem{vonN}
von Neumann, J.:
Discussion on the existence and uniqueness or multiplicity of solutions of the aerodynamical
equation [Reprinted from MR0044302] (1949), {\em Bull. Amer. Math. Soc. {\rm (}N.S.{\rm )}}, {\bf 47} (2010), 145--154.
\end{thebibliography}
\end{document}